%Authors: S. Saccone

%Title: The Pelczynski property for tight subspaces

%Filename: ssacconetight.tex
%TeX: Plain
%Length: 93444 bytes
%Received Date: 12/31/96
%SubjectClass: 46E15
%Abstract: We show that if X is a tight subspace of C(K) then X has
%the Pelczynski property and X^* is weakly sequentially
%complete. We apply this result to the space U of
%uniformly convergent Taylor series on the unit circle
%and using a minimal amount of Fourier theory prove a
%theorem of Bourgain, namely that U has the Pelczynski
%property and U^* is weakly sequentially complete. Using
%separate methods, we prove U and U^* have the
%Dunford-Pettis property. Some results concerning
%pointwise bounded approximation are proved for tight
%uniform algebras. We use tightness and the Pelczynski
%property sto make a remark about inner functions on
%strictly pseudoconvex domains in C^n.

%Citation: Preprint: to appear in Journal Functional Anal.

%Special character check block
%32   space        33 ! exclam. pt.   34 " double quote  35 # sharp
%36 $ dollar       37 % percent       38 & ampersand     39 ' prime
%40 ( left paren.  41 ) rt. paren.    42 * asterisk      43 + plus
%44 , comma        45 - minus 46 . period        47 / divide
%58 : colon        59 ; semi-colon    60 < less than     61 = equal
%62 > greater than 63 ? question mark 64 @ at
%91 [ left bracket 92 \ backslash     93 ] right bracket 94 ^ caret
%95 _ underline    96 ` left single quote
%123 { left brace  124 | vertical bar 125 } right brace  126 ~ tilda

\magnification=\magstep1
\baselineskip=20pt

%
% The Pelczynski Property for Tight Subspaces
%
% Postable Version: get rid of diagrams, use \cal instead of \rsfs.

%
% Scott F. Saccone
% Department of Mathematics
% Washington University
% sfs@math.wustl.edu
%
% Plain TeX Macros
% As of 10/14/96
%
% These are primarily used on my Macintosh.
% Some special inputs:
%
%  (1) the font family "rsfs", a special set of script fonts
%  (2) the macro package "xypic", commutative diagram making software
%
%
%
% Items indigenous to the Mac:
%
% Blackboard Bold font is: msbm10. On my Macintosh, msbm10 is really
% the postscript font "Caslon Openface Regular": the filenames are
% "CasloOpeBTReg" for the postscript outline font, "CasloOpeBTReg.afm"
% for the metric info, and "CasloOpnface BT" for the bitmap screen
% font. I use the EdMetrics tool to accomplish this.
%
% On a unix machine, msbm10 should be the font that comes with
% AMS-Tex fonts.
%
%

%
% Verifications - I put these in various places in my documents where
% there are certain details that may need to be checked over. The
% variable \verifyflag turns these on and off (1=on)
%

\newcount\verifyflag % notes in my papers that I should check something
\verifyflag=0
\def\verify{\ifnum\verifyflag=1\hbox{\bf\ *verify*\ }\fi}

\newcount\commentflag % for comments
\commentflag=0
\def\comment#1{\ifnum\commentflag=1{\bf Comment:\ #1}\fi}

%
% The Kristoffer H. Rose commutative diagram package
%
%
% some possible Xy-pic distribution ftp sites:
%	ftp.mpce.mq.edu.au  in the  /pub/maths/TeX  directory;
%	ftp.diku.dk  		in the  /diku/users/kris/TeX  directory;
%	CTAN in the   tex-archive/macros/generic/diagrams/xypic  directory;
%	  e.g. at  ftp.shsu.edu, ftp.tex.ac.uk  and  ftp.dante.de.
%
% I've created a format file entitled XTtex which takes the place of
% loading this in.
%

% \input xypic

\newcount\showdiagrams
\global\showdiagrams=0 %1=yes 0=no

%
% The following allows the use of Ralph Smith's Formal Script symbols
% in Plain TeX documents.  Use \scr like \cal.
% Set the font sizes and restore the `at' clauses if you want them bigger.
% You can use this method in LaTeX, but only at one basic size.
% If you need symbols in LaTeX titles, captions, etc., work it out or ask
% a LaTeXpert.
%

%
% I'm not sure what a general ftp site would be for the rsfs fonts.
% Right now, rsfs10, rsfs7, and rsfs5 are all "Script12" on my Mac.
%

\def\scr#1{{\cal #1}}

%\font\tenscr=rsfs10 % scaled \magstep1
%\font\sevenscr=rsfs7 % scaled \magstep1
%\font\fivescr=rsfs5 % scaled \magstep1
%\skewchar\tenscr='177 \skewchar\sevenscr='177 \skewchar\fivescr='177
%\newfam\scrfam \textfont\scrfam=\tenscr \scriptfont\scrfam=\sevenscr
%\scriptscriptfont\scrfam=\fivescr
%\def\scr{\fam\scrfam}

%
% Some font definitions
%

\font\csc=cmcsc10          % Caps/small caps
\font\bbb=msbm10           % Blackboard bold (msbm10)
\font\smallroman=cmr10 at 7pt

% Font macros
\def\Bbb#1{\hbox{\bbb #1}}  % Puts things in \bbb (for use in math mode)

%
% Organizational Macros
%

%
% Headlines
%
% Current options: slant text
%

\font\headlinefont=cmcsc10

\def\rightheadline#1{\vbox{
 \line{\hss\sl#1\hss\it\the\count0}\vskip-10pt\line{\hrulefill}}}

\def\leftheadline#1{\vbox{
\line{\hss\sl#1\hss\it\the\count0}\vskip-10pt\line{\hrulefill}}}

\def\monoheadline#1{\vbox{
\line{\hss\headlinefont#1\hss}\vskip-10pt\line{\hrulefill}}}

%
% title
%

\font\titlefont=cmb10 at 15pt
\def\title#1{\centerline{\titlefont#1}\vskip0.5truein}

\font\authorfont=cmcsc10 at 12pt
\def\author#1{\centerline{\authorfont
#1}\vskip0.5truein}

  %comments below the title

%
% abstract
%

\font\smallrm=cmr10 at 8pt
\font\smallcsc=cmcsc10 at 8pt

\font\smallmathit=cmmi10 at 8pt
\font\mathit=cmmi10

\font\scriptrm=cmr10 at 7pt

%
% sections - should be done locally
%

%
% theorem statements, etc.
%
% theorem title in csc, text in italics.
% \endproclaim is the delimiter.
%

\def\proclaimskip{{\vskip10pt}}

\def\endproclaim{\proclaimskip\endgroup}
\def\proclaim#1{\begingroup\parindent=0pt
\proclaimskip
\noindent{\csc#1}\hskip 7pt\it}

%
%  Bibliographic macros
%

\def\semicolon{;}   %to separate the text

\font\bibitfont=cmti10
\font\bibrmfont=cmr10
\font\bibbffont=cmb10

\def\ref#1;#2;#3;#4;#5;#6;{{\bibrmfont \item{[#1]} #2, {\bibitfont #3,}
#4 {\bibbffont #5} #6.}}
     %number,author(s),title,journal,volume,whatever.

\def\book#1;#2;#3;#4;{{\bibrmfont \item{[#1]} #2, {\bibitfont #3,} #4.}}
     %number,author(s),title,whatever.
\def\bysame{\hbox to 30pt{\hrulefill}}

%
% Text macros - General
%

%
% Teaching Macros
%

% Example:
% \points{15}a) Let $\u=(1,2,3)$ and $\v=(3,2,1)$. Find the equation...
%

%
%  Text macros - Mathematical
%

% General Math
\def\bbox{\vrule height5pt width5pt depth0pt} %End of proof box.

\def\claim{{\par\csc Claim:\ }}
\def\claimII{{\csc Claim:\ }}

\def\iff{if and only if}
\def\implies{{$\,\Rightarrow\,$}}

\def\proof{\vskip-\parskip{\it Proof\ \ }}
\def\pclaim#1{{\csc\par Claim #1:\ }}
  % To be used in a sentence
\def\qed{\ \hfill\bbox\vskip\baselineskip}
\def\tfae{the following are equivalent:}
\def\varproof#1{{\vskip10pt{\it #1.}}}
\def\Wlog{{without loss of generality}} % I think \wlog is taken
\def\wrt{with respect to}

% General Symbols
\def\iffarrow{{$\,\Leftrightarrow\,$}}

% Mathematical Names
\def\leb{Lebesgue}
\def\pel{Pe\l czy\'nski}

% Banach Space Theory
\def\bs{Banach space}
\def\ccont{completely continuous}

\def\cseq{$c_0$-sequence}
\def\dpp{Dunford-Pettis property}
\def\fr{{finitely representable}}

\def\lonesp{{$\lone$-space}}
\def\loneseq{$\llone$-sequence}

\def\pp{\pel\ property}
\def\pwae{{pointwise~a.e.\thinspace}}
\def\rwc{relatively weakly compact}
\def\sepdist{separable distortion}
\def\sepdistone{\sepdist\ of an $L^1$-space}

\def\ucc{unconditionally converging}

\def\wsc{weakly sequentially complete}

\def\wuc{weakly unconditionally Cauchy}

% Uniform Algebras

\def\ntmrb{nontrivial minimal reducing band}

%
%  Mathematical macros
%

%General
\def\biggoesto{{\,\longrightarrow\,}}
\def\bigmod#1{{\bigl|#1\bigr|}}
\def\bigrestrict#1;#2{#1\kern-.5pt\lower1.8pt\hbox
 {$\big |_{\lower1pt\hbox{$\scriptstyle\kern-1pt#2$}}$}}

\def\brace#1{{\left\{ #1 \right\}}}
\def\Cn{{\C^n}}
\def\C{{\Bbb{C}}}

\def\eps{{\varepsilon}}
\def\exp#1{\hbox{e}^{#1}}
\def\goes#1{{\,\buildrel #1\over\longrightarrow\,}}
\def\goesto{{\,\longrightarrow\,}}
\def\iint{\mathop{\int\kern-6pt\int}}
\def\iiint{\mathop{\int\kern-4pt\int\kern-4pt\int}}

\def\inverse#1{{#1^{-1}}}
\def\iprod#1{\langle #1 \rangle}

\def\modulus#1{{\left|#1\right|}}
\def\mpoly{{\Bbb B^m}}
\def\nball{{\Bbb B_n}}
\def\notgoesto{{\goesto\kern-12pt/\kern7pt}}
\def\parenth#1{{\left(#1\right)}}
\def\pint#1;#2{{\int #1\,d#2}}

\def\pseq#1#2{{\brace{#1_#2}}}
\def\R{{\Bbb R}}
\def\restrict#1;#2{#1\kern-.5pt\lower.5pt\hbox{$|_{\lower1pt\hbox
 {$\scriptstyle\kern-1pt#2$}}$}}
\def\seq#1{{\brace{#1_n}}}

\def\text#1{{\quad\hbox{#1}\quad}}
% This makes a limsup symbol with the bar over a sup. The vphantom
% control sequence tells tex to put the symbol in a box who
% depth is the same as that of `sup' (the `p' makes it lower).
% underneath, it is not too high.
\def\varlimsup{\mathop{\vphantom{\rm sup}\overline{\rm lim}}}

% Calculus and Vectors

\def\tricolvect#1;#2;#3;{\pmatrix{#1\cr#2\cr#3}}
\def\quadcolvect#1;#2;#3;#4;{\pmatrix{#1\cr#2\cr#3\cr#4}}

%Complex Variables
\def\bidisk{{\overline\Delta\times\overline\Delta}}
\def\dbar{\bar\partial}

%Functional Analysis
\def\bigllplus{\mathop{\,\bigoplus\nolimits_\llone\,}}
\def\bmd{{\hbox{d}}} % Banach-Mazur distance
\def\chull{{\hbox{\rm co}}}
\def\clspan{{\overline{\rm sp}\,}}
\def\ddual#1{{#1^{**}}}
\def\dist{{\hbox{\rm dist}}}
\def\goesstar{{\goes{w^*}}}
\def\goesweak{{\goes w}}
\def\hulli#1#2{{\chull\brace{#1_1,#1_2,\ldots,#1_#2}}}
\def\hullii#1#2{{\chull\brace{#1_{#2+1},#1_{#2+2},\ldots}}}

\def\injects{{\,\hookrightarrow\,}}

\def\linf{{L^\infty}}
\def\llinf{{l^\infty}}

\def\llone{{l^1}}
\def\llplus{{\,\oplus_\llone\,}}

\def\lone{{L^1}}
\def\lonedirect{{\sum_{l^1}\kern-1.5pt\oplus}}

\def\net#1{{\brace{#1_\alpha}}}
\def\norm#1{\Vert #1 \Vert}
\def\normcl#1{{\overline{#1}^{\hbox{\smallroman norm}}}}

\def\txy{{T:X\to Y}}
\def\wkstarcl#1{{\overline{#1}^{w^*}}}

% Related to the Pelczynski property

\def\wucprodii#1#2{{#1_{n}\prod\limits_{#2=1}^{n-1}\bigl(1-|#1_#2|)}}

%Uniform Algebras
\def\aperp{{A^\perp}}
\def\baperp{{\bd_{\kern-1.8pt\aperp}}} %band gen by aperp.
\def\bxperp{{\bd_{\kern-1.8pt\xperp}}}
\def\bd{{\scr B}}
\def\bmodaperp{{\baperp\kern-1pt/\kern-1pt\aperp}}
\def\bmodxperp{{\bxperp\kern-1pt/\kern-1pt\xperp}}
\def\hinf{{H^\infty}}

\def\maxidspof#1{{{\scr M}_{\kern-1pt\lower2pt\hbox{$\scriptstyle{#1}$}}}}
\def\sb{{\scr S}}
\def\shilov#1{{\partial_{\kern-1.5pt#1}}}
\def\xperp{{X^\perp}}
%\def\varbaperp{{\bd\over\bd\cap\aperp}}

%Measure Theory
\def\ldecomp#1{{{#1}_a+{#1}_s}}

%
% end of my main macro file
%

% \input xypic

\showdiagrams=0

\verifyflag=0  % for the special "verify" symbols I use to check details
\commentflag=0

\headline{\ifnum\count0>1\monoheadline{Tight Subspaces}\fi}

%
% special macros for this paper
%

% math

\def\MplusL{{M\llplus\kern-1.5pt\lone(\mu)}}
\def\Nth{N^{\hbox{\smallroman th}}}

\def\Svarlaperptoba#1{{\varlaperp{#1}\,{\buildrel\sigma\over\longrightarrow}
 \,\varbmodaperp}}
\def\Thbatoh{{\hinf(\baperp)\,{\buildrel\tau\over\longrightarrow}\,\hinf}}

\def\XmodE{{X\over E}}

\def\YmodYE{{Y\over Y\cap E}}
\def\YtoE{{\YmodYE\,\buildrel\sigma\over\longrightarrow\,\XmodE}}

\def\apb#1#2{{#1\kern-1pt+#2}}
\def\barf{\bar f}
\def\bmodbaperp{{\bd/\bd\cap\aperp}}
\def\boundaryD{{\partial D}}
\def\falpha{f_\alpha}

\def\lamq{{\lambda_Q}}

\def\mpl{{M\llplus L}}
\def\mrbsum{{\bigllplus\bd_\alpha/\bd_\alpha\cap\aperp}}
\def\ph{\varphi}
\def\sa{{\sb_{\aperp}}}
\def\subnot{{\scriptscriptstyle 0}}
\def\varBaplusSa{{\varbmodaperp\llplus\sa}}

\def\varbmodaperp{{\baperp\over A^\perp}}

\def\varlaperp#1{{L^1(#1)\over L^1(#1)\cap\aperp}}

\def\vf{{\varphi}}

% organizational
\def\section#1{\vskip1.5\baselineskip\centerline{\csc
#1}\vskip1.5\baselineskip}

%
% theorems
%

% section 1
%\def\resIa{Lemma 1.a}          % on wuc series -> c_0 seq
%\def\resIb{Proposition 1.b}    % Chaumat's result
%\def\resIc{Theorem 1.2}        % Kisliakov's theorem
\def\resId{Proposition 1.1}    % tfae for pp

% section 2
\def\resIIa{Theorem 2.1}       % thm tight>pp
\def\resIIb{Theorem 2.2}       % thm R.C. James
\def\resIIc{Lemma 2.3}         % temma 1
\def\resIId{Lemma 2.4}         % lemma 2

% section 3
\def\resIIIa{Lemma 3.1}        % on reducing bands/unions S_g
\def\resIIIb{Proposition 3.2}  % tight>b/a sep
\def\resIIIc{Theorem 3.4}      % sep dists
\def\resIIId{Lemma 3.5}        % on H^\inf(bd) and self-adj
\def\resIIIe{Theorem 3.7}      % Lind-Pel
\def\resIIIf{Proposition 3.8}  % products->sep-dists
\def\resIIIg{Proposition 3.6}  % \bmod not embed in L^1
\def\resIIIh{Lemma 3.10}       % lemma on approx in hinf(bd)
\def\resIIIi{Proposition 3.11} % prop. davie<>quo. map
\def\resIIIj{Proposition 3.12} % equiv on weak rich
   % def full (not used)
\def\resIIIl{Proposition 3.13} % contruct full subsp
\def\resIIIm{Proposition 3.14} % full bands
\def\resIIIn{Corollary 3.15}   % to prop 3.m i
\def\resIIIo{Corollary 3.16}   % ii
\def\resIIIp{Theorem 3.18}     % tight>str Davie
\def\resIIIq{Corollary 3.17}   % iv
\def\resIIIr{Corollary 3.3}    % corr to prop 3.2
      % Kisliakov's theorem (not used)
\def\resIIIt{Corollary 3.9}    % products->s.d.
\def\resIIIv{Theorem 3.19}     % Mooney extensions-A(D)

%section 4 - added after a suggestion by the referee
\def\resIVna{Theorem 4.1}      % all about U
\def\resIVnb{Proposition 4.2}  % X is tight
\def\resIVnc{Lemma 4.3}        % S^** is cc
\def\resIVnd{Proposition 4.4}  % X_B=C
\def\resIVne{Proposition 4.5}  % X_CG=Y

% section 5 - was originally section 4
\def\resIVa{Theorem 5.1}       % thm on innter fncs
\def\resIVb{Theorem 5.2}       % lemma tau isometry
\def\resIVc{Lemma 5.3}         % lemma

%
% equations
%

% section 2
\def\eqnIIa{(2.1)}   % in lemma IId

% section 3
\def\eqnIIIa{(3.1)}  % in lemma IIIa
\def\eqnIIIb{(3.2)}  % in lemma IIId
\def\eqnIIIc{(3.3)}  % in sdprop
\def\eqnIIId{(3.4)}  % in sdprop
\def\eqnIIIe{(3.5)}  % ""
\def\eqnIIIf{(3.6)}  % ""
\def\eqnIIIg{(3.7)}  % ""
\def\eqnIIIh{(3.8)}  % ""
\def\eqnIIIi{(3.9)}  % ""
\def\eqnIIIj{(3.10)} % ""
\def\eqnIIIk{(3.11)} % ""
\def\eqnIIIl{(3.12)} % hinf(b)->hinf(m)
\def\eqnIIIm{(3.13)} % predual of above
\def\eqnIIIn{(3.14)} % def of full
\def\eqnIIIo{(3.15)} % def of full ii

% section 4
\def\eqnIVna{(4.1)}  % V is 1-summing
%
% references
%

% removed:
%\def\refBD{[4]}
%\def\refCCG{[5]}
%\def\refDS{[12]}
%\def\refPrieto{[16]}

\def\refBP{[1]} % Bessaga Pelczynski
\def\refBourU{[2]}
\def\refBourPP{[3]}
\def\refBourDPP{[4]}
\def\refCT{[5]}   % Cima and Timoney
\def\refCG{[6]}
\def\refCR{[7]}
\def\refCon{[8]}
\def\refDav{[9]}
\def\refDeli{[10]}  % Delbaen: disk algebra
\def\refDelii{[11]} % DelbaenL R(K)
\def\refGamUA{[12]}
\def\refHeni{[13]}
\def\refJames{[14]}
\def\refKisPP{[15]}
\def\refKisUA{[16]}
\def\refLP{{[17]}}
\def\refMilne{[18]}
\def\refMoon{[19]}
\def\refPP{[20]} % Pelczynski: the PP
\def\refPeli{[21]}  % BS of analytic
\def\refMe{[22]}
\def\refWojPP{[23]}
\def\refWojBS{[24]}

%
% \rsfs substitution
%

\def\scr#1{{\cal #1}}

\title{The \pel\ Property for Tight Subspaces}
\author{Scott F. Saccone}
\vskip20pt

\global\textfont0=\smallrm
\global\scriptfont0=\scriptscriptfont0
\global\textfont1=\smallmathit

\centerline{\vbox{\hsize=5.5truein\baselineskip=10pt{\smallcsc Abstract.}
\smallrm
We show that if $X$ is a tight subspace of $C(K)$ then $X$ has the \pp\ and
$X^*$
is weakly sequentially complete. We apply this result to the space $U$
of uniformly convergent Taylor series on the unit circle and using a minimal
amount of Fourier theory prove that $U$ has the \pp\ and $U^*$ is weakly
sequentially
complete. Using separate methods, we prove $U$ and $U^*$ have the \dpp.
Some results
concerning pointwise bounded approximation are proved for tight uniform
algebras. We
use tightness and the \pp\ to make a remark about inner functions on strictly
pseudoconvex domains in
$\Cn$.}}
\vskip30pt

\global\textfont0=\tenrm
\global\scriptfont0=\scriptrm
\global\textfont1=\mathit

\section{1. Introduction and Background}

The \pp, whose concept was influenced by the
work of Orlicz, was introduced by \pel\ in \refPP. We say a sequence $\seq
x$ in a
\bs\ $X$ is a {\wuc\ } series (w.u.C. series) if $\sum|x^*(x_n)|<\infty$
for every
$x^*\in X^*$ and $\seq x$ is an {\ucc\ series} if $\sum x_{\pi(n)}$
converges in
norm for every permutation $\pi$ of the natural numbers. If $X$ and $Y$ are
\bs s and
$T:X\to Y$ is a continuous linear operator we say $T$ is an {\it\ucc\ }
operator
if $T$ takes every \wuc\ series in $X$ to an \ucc\ series in $Y$. It follows
from the work of Orlicz that every weakly compact operator is an \ucc\
operator.

The \pp\ for a \bs\ is the realization of a converse to the result of
Orlicz. We
say $X$ has the \pp\ if every \ucc\ operator on $X$ is weakly compact. It is a
theorem of Bessaga and \pel\ in \refBP\ \verify\ that a continuous linear
operator $\txy$ is \ucc\ \iff\ $T$ is never an isomorphism on a copy of $c_0$
in $X$.

We say a sequence $\seq x$ is a \cseq\ if it is a basic sequence which is
equivalent to the unit vector basis of $c_0$ and similarly for \loneseq s.
Given a bounded subset $E\subset X^*$ we will be interested in knowing when
there exists a \wuc\ series $\seq x$ in $X$ that fails to tend to zero
uniformly on $E$; that is, $\varlimsup\limits_{n\to\infty}\sup\limits_{x^*\in
E}|x^*(x_n)|>0$.  It follows from the result of Bessaga and \pel\ mentioned
above that this is equivalent to the existence of a $c_0$ sequence $\seq x$
that fails to tend to zero uniformly on $E$ (just consider the operator
$T:X\to\llinf(E)$ by $Tx(x^*)=x^*(x)$). We say a sequence $\seq x$ is a
weak-Cauchy
sequence if $\lim x^*(x_n)$ exists for every $x^*\in X^*$ and we say a \bs\
$X$ is {\it\wsc} if every weak-Cauchy sequence in $X$ is weakly convergent.

The following are some more or less well-known characterizations of the \pp.

\proclaim{\resId} If $X$ is a \bs\ then the following are equivalent.

(a) $X$ has the \pp.

(b) If $\txy$ is a continuous linear operator which fails to be weakly compact
then $T$ is an isomorphism on some copy of $c_0$ in $X$.

(c) If $E\subseteq X^*$ and the weak closure of $E$ fails to be weakly
compact then
there exists a \wuc\ series $\seq x$ in $X$ which fails to tend to zero
uniformly on
$E$.

(d) The following hold: (i) $X^*$ is \wsc\ (ii) If $\seq{x^*}$ is an \loneseq\
in $X^*$ then there exists a \cseq\ in $X$ such that
$\modulus{x^*_{n_k}(x_k)}>\delta>0$ for all $k$ for some sequence
$\brace{n_k}$.
\endproclaim

The equivalence of (a) and (b) follows from the theorem of Bessaga and \pel\
mentioned above, while the equivalence of (a) and (c) is well-known. The
equivalence of (a) and (d) is less popular, but can be deduced from
(c) and the now ubiquitous result of Rosenthal and Dor: if $X$ is any \bs\ and
$\seq x$ is a bounded sequence in $X$ which has no weak-Cauchy subsequence
then $\seq x$ has an $\llone$-subsequence.

\comment{check that the above observations aren't trivial and verify the
properties of the James space}

All $C(K)$ spaces were shown to have the \pp\ in \refPP.
Every infinite-dimensional \lonesp\ fails to have this property since these
spaces do not contain a copy of $c_0$. Delbaen and Kisliakov independently
showed the disk algebra has the \pp\ in \refDeli\ and \refKisPP\ respectively.
Delbaen extended these results to~$R(K)$ for special classes of planar sets~$K$
in \refDelii\ as did Wojtaszczyk in \refWojPP. It was shown that $R(K)$ has the
\pp\ for every compact planar set $K$ in \refMe. It was also shown in \refMe\
that every so-called T-invariant uniform algebra on a compact planar set has
the \pp. The T-invariant class includes $R(K)$ as well as $A(K)$ for all
compact
planar sets $K$. However, it is not known if any of these planar
uniform algebras fail to be linearly isomorphic to the disk algebra. Bourgain
showed the ball-algebras and the polydisk-algebras have the \pp\ in \refBourPP.
This result was extended in \refMe\ to $A(D)$ for strictly pseudoconvex domains
$D$ in $\Cn$.

Not all uniform algebras have the \pp. In fact, it is a result of Milne in
\refMilne\ that every \bs\ $X$ is isomorphic to a complemented subspace of a
uniform algebra $A$, where $A$ can be taken to be the uniform algebra
on $B_{X^*}$ generated by $X$. However, the author is not presently aware of
any uniform algebras on compact subspaces of $\R^n$ which fail to have the \pp.

The \pp\ holds for a special class of spaces which includes many
examples of uniform algebras of analytic functions. If $K$ is a compact
space and
$X\subseteq C(K)$ is a closed subspace then we say $X$ is a {\it tight
subspace} if
the operator $S_g:X\to C(K)/X$ by $f\mapsto fg+X$ is weakly compact for
every $g\in
C(K)$. We say a uniform algebra $A$ on $K$ is a {\it tight uniform algebra}
if it is
a tight subspace. The concept of tightness was introduced by B. Cole and T.W.
Gamelin in \refCG\ as the ability to solve an abstract $\dbar$-problem with a
mild gain in smoothness.

Although the authors in \refCG\ were mainly interested in weakly compact
Hankel-type operators, in many of the examples the operators
$S_g$ were proven to be compact. We say $X\subseteq C(K)$ is a {\it strongly
tight subspace} if $S_g$ is compact for every $g$, and similarly we define
strongly tight uniform algebras. It was proven in \refCG\ that $R(K)$ is
strongly tight for every compact planar set $K$, and also $A(D)$ is strongly
tight for every strictly pseudoconvex domain $D$ in $\Cn$ with $C^2$ boundary.
More generally, $A(D)$ will be strongly tight whenever the $\dbar$ problem can
be solved in $D$ with H\"older estimates on the solutions.

Currently there is no known example of a tight uniform algebra which fails to
be strongly tight. However there are examples of tight, non-strongly
tight subspaces.  We say an operator
$\txy$ between \bs s is {\it\ccont\ } is $T$ takes weakly null sequences to
norm
null sequences. We say a \bs\ $X$ has the {\it\dpp} if every weakly compact
operator $\txy$ is completely continuous. It follows from the work of
Bourgain in
\refBourDPP\ that any strongly tight subspace has the \dpp. By considering the
identity operator, we see that every infinite-dimensional reflexive space fails
to have the \dpp. Hence, any infinte-dimensional reflexive space $X$ will
be tight
in any $C(K)$ space it is embedded in, but can never be realized as a strongly
tight subspace.

Our main result is \resIIa\ which states that every tight
subspace of $C(K)$ has the \pp. This result generalizes a theorem from \refMe,
namely that every strongly tight uniform algebra has the \pp\ (actually,
the proof
in \refMe\ does not use any algebraic structure and would work for strongly
tight
subspaces of $C(K)$).

An application of this theorem was kindly forwarded to the author by the
referee. Let $U$ be the space of continuous functions on the unit circle which
extend to be analytic in the unit disk and whose Taylor series converge
uniformly on
the closed disk. We define a norm on $U$ by taking the supremum of
sup-norms the
partial sums of the Taylor series. Included in the referee's report was a
fairly short
proof that $U$ embeds into some $C(K)$ space as a tight subspace. Using a
shortcut that
allows us to check that $S_g$ is weakly compact for a only small collection of
functions $g$, we give an even simpler proof of this result in \resIVnb. It
now follows
from \resIIa\ that $U$ has the \pp. This is a result that has been
established by
Bourgain in \refBourU. Bourgain's proof uses a fair amount of hard
analysis, including
Carleson's theorem on the pointwise almost everywhere convergence of
Fourier series
in $L^2$. Interestingly enough, our proof uses little more than the
Plancherel theorem.

The results on the space $U$ are in Section 4. In addition to proving $U$
has the \pp,
we show that $U$ and $U^*$ have the \dpp. The main ingredient is a
theorem of Bourgain which concerns the operators $S_g$. Bourgain proves in
\refBourDPP\ that a closed subspace $X\subseteq C(K)$ will have a dual
space with
the \dpp\ whenever $\ddual{S_g}$ is completely continuous for every $g\in
C(K)$. It is
well-known that a \bs\ $Y$ has the \dpp\ whenever $Y^*$ does. We show that
$U$ embeds into a $C(K)$ space (the same $K$ indicated above) as a subspace $X$
satisfying Bourgain's criteria. As in the tight subspace case, the proof is
quite
simple, and uses very little Fourier theory. Our results on the \bs\
structure of
$U$ are summarized in \resIVna. Using the same proof used for the disk algebra
(see \refPeli), we prove the known result that $U$ is not isomorphic to a
quotient
of $C(G)$ for any compact space $G$.

We noted above that $U$ is isomorphic to a tight subspace $X$ of $C(K)$ for
some space
$K$. We prove that $X$ is not strongly tight and give a characterization of
those $g$
for which $S_g$ is compact. Hence, in addition to the reflexive spaces, $X$
yields a
new example of a tight, non-strongly tight space. As we noted above, an
example of a
tight algebra of functions which is not strongly tight has not yet been
produced.

In addition to the \pp, we investigate properties of tightness which are more
commonly studied in the context of function algebras. Let $A$ be a uniform
algebra and let $\maxidspof A$ be the maximal ideal space of $A$. If
$\varphi$ and $\psi$ are elements of $\maxidspof A$ then we say $\varphi$ and
$\psi$ are in the same {\it Gleason part} if $\norm{\varphi-\psi}_{A^*}<2$.
This is an equivalence relation where the classes are called the Gleason parts
of $A$. We say a part is trivial if it consists of one point. It was shown in
\refCG\ that every tight uniform algebra on a compact metric space $K$
possesses at most countably many nontrivial Gleason parts. We give a simple
proof of this fact. We will need the theory of bands of measures (for more
information on bands and related ideas see \refCG\ or \refCon).

Let~$K$ be a compact space. If~$\bd\subseteq M(K)$ we
say $\bd$ is a {\it band of measures} if~$\bd$ is a closed subspace of~$M(K)$
and when~$\mu\in\bd$,~$\nu\in M(K),$~and $\nu\ll\mu$
then $\nu\in\bd.$ The Lebesgue Decomposition Theorem says that if $\mu\in M(K)$
then $\mu$ can be uniquely written as $\mu=\mu_a+\mu_s$ where $\mu_a\in\bd$
and $\mu_s$ is singular to every element of $\bd.$ If $\bd$ is a band, the
{\it complementary band} $\bd^\prime$ of $\bd$ is the collection of measures
singular to every measure in $\bd.$ It follows that $M(K)=\bd\llplus\bd'.$
It is
a well-known fact that if $\bd$ is a band then $\bd\cong\lone(\mu)$ for some
abstract measure $\mu$.

If $\bd$ is a band, we define $\linf(\bd)$ to be the space of uniformly bounded
families of functions $F=\{F_\nu\}_{\nu\in\bd}$ where $F_{\nu}\in\linf(\nu)$
and $F_\nu=F_\mu$ a.e. $[d\nu]$ whenever $\nu\ll\mu$.  The norm in $\linf(\bd)$
is given by $\norm F=\sup_{\nu\in\bd}\Vert F_\nu\Vert_{\linf(\nu)}.$
The pairing $\left <\nu,F\right >=\int F_\nu\,d\nu$ defines an isometric
isomorphism between $\linf(\bd)$ and $\bd^*.$ If $X$ is
a subspace of $C(K)$ let $\hinf(\bd)$ and $\hinf(\mu)$ be the
weak-star closure of $X$ in $\linf(\bd)$ and $\linf(\mu)$ respectively.
If~$\mu\in\bd$, there is a natural projection
$H^\infty(\bd)\to H^\infty(\mu)$ defined by $F\mapsto F_\mu$. We define
$\bxperp$ to be
the band generated by the measures in $\xperp$ and $\sb$ to be the band
complement to
$\bxperp$. It follows from the Lebesgue decomposition that
$X^*\cong\bmodxperp\llplus\sb$.

We say a band~$\bd$ is a {\it reducing band} for~$X$ if for any
measure~$\nu\in\xperp$ the projection~$\nu_a$ of~$\nu$ into~$\bd$ by
the \leb\ decomposition is also in~$\xperp.$ We say~$\bd$ is a {\it minimal
reducing band} if~$\bd\ne\brace0$ while~$\brace{0}$ is the only
reducing band properly contained in~$\bd.$

Suppose, for now, the subspace $X$ is a uniform algebra $A$. The following
version of
the Abstract F. and M. Riesz Theorem can be found in \refCG.
Let~$\varphi\in\maxidspof
A.$ Then the band generated by the representing measures for~$\varphi$ is a
minimal
reducing band. The band generated by the representing measures
for~$\varphi$ is equal
to the band generated by the representing measures for all the points in
the same
Gleason part as $\varphi$. Hence every Gleason part of a uniform algebra
corresponds to
a distinct minimal reducing band.

If $A$ is a uniform algebra we say a point~$z\in K$ is a {\it peak point}
for~$A$ if
there exists an element~$f\in A$ such that~$f(z)=1$ and~$|f(w)|<1$
for~$w\ne z.$ We
say~$z$ is a {\it generalized peak point} if the only complex representing
measure for~$z$ is the point mass at~$z.$ The {\it Choquet boundary} of~$A$ is
the collection of all generalized peak points. The point-evaluations for the
points off the Choquet boundary lie in~$\bmodaperp$ while those for the points
on the Choquet boundary lie in~$\sb.$

If $\bd$ is a minimal reducing band for $A$ and $\bd\subseteq\sb$ then it
can be seen
that $\bd$ is all multiples of a point mass~$\delta_z$ at some generalized
peak point~$z\in K.$ We call these {\it trivial minimal reducing bands} and the
others {\it nontrivial minimal reducing bands.} Note that a minimal reducing
band~$\bd$ is trivial \iff\ $\bd\cap\aperp=0$ (this implies every subband
of $\bd$
is a reducing band).

Note that the intersection of two reducing bands is a reducing band and so two
minimal reducing bands either coincide or are singular. If we let~$\net\bd$
be the
collection of all the non-trivial minimal reducing bands for $A$ then
$\bigllplus\bd_\alpha$ is a reducing band contained in~$\baperp.$ However,
this may
not be all of~$\baperp.$ For more information, see \refCG. The sum
$\bigllplus\bd_\alpha/\bd_\alpha\cap\aperp$ is now isometric to a closed
subspace
of~$A^*$ which is contained in~$\bmodaperp.$

We show in Section 3 that if $X$ is a tight subspace of $C(K)$ for a
metric space $K$ then $\bmodxperp$ is separable. It then follows that $X^*$
is a
separable distortion of an $\lone$-space; that is, the dual of $X$ can be
written as
the direct sum of an $\lone$-space and a separable space. When $A$ is a
tight uniform
algebra on a metric space $K$ then we see that $A$ can have at most
countably many
nontrivial Gleason parts and at most countably many \ntmrb s. This
conclusion is
easily deduced once we see that $\bmodaperp$ is separable. In fact we show
that $A$
will have at most countably many nontrivial Gleason parts whenever $A^*$ is
merely
embedded in a separable distortion of an $\lone$-space. The proof is an
adaptation on
a method of Henkin which can be found in \refPeli.

The separability of $\bmodaperp$ has some interesting consequences. For
example, when
$K$ is a metric space this implies that there exists an $m\in\baperp$ such
that every
non-peak point of $A$ has a representing measure absolutely continuous
\wrt\ $m$.
Furthermore, $m$ will have some other special properties concerning pointwise
bounded approximation. The prototypical example is the following. Let $K$
be a compact
subspace of $\C$. Let $Q\subseteq K$ be the non-peak points of $R(K)$ and let
$\lambda_Q$ be \leb\ planar measure restricted to $Q$. Let
$\hinf(\lambda_Q)$ be the
weak-star closure of $R(K)$ in $\linf(\lambda_Q)$. It is a theorem of A.M.
Davie
in \refDav\ that given $f\in\hinf(\lambda_Q)$ there exists a sequence
$f_n\in R(K)$
with $\norm{f_n}\leq\norm f$ such that $f_n\goesto f$ pointwise a.e.
$[\lambda_Q]$.
The measure $m$ will possess the same property as $\lambda_Q$.

Section 5 uses the theory of tight uniform algebras and the \pp\ to deduce
a result about inner functions on strictly pseudoconvex domains in $\Cn$. The
background is as follows. We say a subalgebra $B\subseteq\linf(m)$, where
$m$ is \leb\
measure on the unit circle, is a Douglas algebra if $B$ contains $\hinf$.
Recall the
Chang-Marshall Theorem which states that every Douglas algebra on the unit
circle is
generated by $\hinf$ and a collection of conjugates of inner functions. In
contrast
to this result we prove the following. If $n>1$ and then there are no
nonconstant
inner functions whose conjugate is in $\hinf(m)+C$ where $m$ is surface
area measure
on the unit ball in $\Cn$. It is well-known that $\hinf(m)+C$ is a closed
subalgebra
of $\linf(m)$. The proof is identical when the unit ball is replaced by a
strictly
pseudoconvex domain $D$ which has $C^2$ boundary. The proof is quite soft
and relies
mainly on the fact that $A(D)$ is a strongly tight subalgebra of $C(\bar D)$.

\section{2. The \pp}

In this section we prove our main result concerning tight subspaces of $C(K)$.

\proclaim{\resIIa} Let $X\subseteq C(K)$ be a tight subspace. Then $X$ has
the \pp.
\endproclaim

The following well-known theorem on weak compactness will be essential.
Recall that
a set $E$ is {\it relatively weakly compact} if the weak closure of $E$ is
weakly
compact.

\proclaim{\resIIb\ (R.C. James)} Let $X$ be a \bs\ and
let $E\subset X$ be a bounded subset. Then \tfae

(a) $E$ fails to be \rwc.

(b) There exists a sequence $\seq x$ in $E$ and a $\rho>0$ such that
if $$\eqalignno{V_n&=\hulli xn\cr
\noalign{and}
W_n&=\hullii xn}$$ then $\dist(V_n,W_n)>\rho$ for
all $n.$

(c) There exist sequences $\seq\varphi$ in $B_{X^*}$ and $\seq x$ in $E$
and a $\rho>0$ such that
$$\eqalign{
\varphi_n(x_k)&=0\quad\hbox{\rm for\ }1\leq k\leq n\cr
\hbox{\rm Re}\,\varphi_n(x_k)&>\rho\quad\hbox{\rm for\ }k\geq n+1.}
$$
\endproclaim

The following lemmas deal with non-weakly compact sets in arbitrary \bs s.
The second
lemma is an integral part of the gliding hump construction used to prove
\hbox{\resIIa}.

\proclaim{\resIIc} If $\txy$ is a continuous linear operator
and $S:X\to Z$ is weakly compact and $\ddual x\in B_{\ddual X}$
with $\norm{\ddual T\ddual x}>\rho>0$ and $\norm{\ddual S\ddual x}<\eps$
then there exists an $x\in B_X$ with $\norm{Tx}>\rho$ and $\norm{Sx}<\eps.$
\endproclaim

\proof Choose $z^*\in Z^*$ with $\norm{z^*}=1$
and $\iprod{\ddual T\ddual x,z^*}>\rho$. Let $x^*=T^*z^*$ and
define $\Omega=\brace{x\in B_X\,\bigm|\,\hbox{Re}\,\iprod{x,x^*}>\rho}$.
Then $\Omega$ is convex with $\ddual x\in\wkstarcl\Omega$
and $\norm{Tx}>\rho$ for all $x\in\Omega.$ Since $S$ is weakly
compact we have $\wkstarcl{S(\Omega)}=\normcl{S(\Omega)}$ and
so $\ddual S\ddual x\in\normcl{S(\Omega)}$.
Since $\norm{\ddual S\ddual x}<\eps$ we may find $x\in\Omega$
with $\norm{Sx}<\eps.$\qed

\proclaim{\resIId} Let $X$ be a \bs\ and suppose $E\subset X^*$ is bounded
and fails to be relatively weakly compact. Then there exists a $\rho>0$
and a subset $F\subseteq E$ such that for any infinite
subset $F'\subseteq F$ and any weakly compact linear
operator $S:X\to Z$ there exist sequences $\seq x\subset B_X$ and
$\seq{\zeta^*}\subseteq F'$ with $\bigmod{\zeta_n^*(x_n)}>\rho$ and
$\norm{Sx_n}\goesto0.$
\endproclaim

\begingroup % proof
\def\ph{\varphi}

\proof Assume $E\subset X^*$ is bounded and fails to be relatively
weakly compact. Then by the R.C. James theorem there exists a $\rho'>0$
and a sequence $\seq{x^*}\subseteq E$ and $\seq\ph\subset B_{\ddual X}$
so that
$$\eqalign{
\ph_n(x^*_k)=0\quad\hbox{for\ }1\leq k\leq n\cr
\hbox{Re}\,\ph_n(x^*_k)>\rho'\quad\hbox{for\ }k\geq n+1.}\eqno\eqnIIa$$
Let $F=\seq{x^*}.$ Suppose $F'\subset F$ is infinite. Without loss
of generality we may assume $F'=F;$ that is, $F'$ will be a sequence
satisfying \eqnIIa\ with a subsequence of the $\seq\ph$ but with the same
constant $\rho'.$

Let $Y=\llinf(F')$ and let $T:X\to Y$ be the canonical map.
If $x^*\in F'$ define~\hbox{$\delta_{x^*}\in\llone(F')$} to be the point mass
at $x^*$ so $(\ddual T\ph_n)(\delta_{x_k^*})=\ph_n(x^*_k)$.
Let $\ddual{y_n}=\ddual T\ph_n$ and let $y^*_n=\delta_{x^*_n}.$ Then
$$\eqalign{
\iprod{y^*_k,\ddual{y_n}}&=0\quad\hbox{for\ }n\geq k\cr
\hbox{Re}\,\iprod{y^*_k,\ddual{y_n}}&>\rho'\quad\hbox{for\ }1\leq n\leq k-1}
$$
for $n\geq1.$ Therefore
$$\dist\Bigl(\chull\brace{\ddual{y_1},\ddual{y_2},\ldots,\ddual{y_{k-1}}},
\chull\brace{\ddual{y_k},\ddual{y_{k+1}},\ldots}\Bigr)>\rho'$$ for $k\geq2.$

Choose $$\ddual{u_n}\in\chull\brace{\ph_1,\ph_2,\ldots,\ph_n}$$
and $$\ddual{v_n}\in\chull\brace{\ph_{n+1},\ph_{n+2},\ldots}$$ so
that $\norm{\ddual T(\ddual{u_n}-\ddual{v_n})}>\rho'.$ Since $\ddual S$
is weakly compact we have
$\norm{\ddual S(\ddual{u_n}-\ddual{v_n})}\goesto0.$
Let $\ddual{x_n}={1\over2}(\ddual{u_n}-\ddual{v_n})$
so $\ddual{x_n}\in B_X$ and $\norm{\ddual S\ddual{x_n}}\goesto0$
and $\norm{\ddual T(\ddual{x_n})}>\rho'/2.$ Let $\rho=\rho'/2$ and
choose $x_n\in B_X$ by~\resIIc\ so that $\norm{Sx_n}\goesto0$
and $\norm{Tx_n}>\rho.$ By definition of $T$ we may find $\zeta^*_n\in F'$ so
that $\bigmod{\zeta_n^*(x_n)}>\rho.$ This completes the proof.\qed

\endgroup % proof

We now return to the \pp.

\vskip\baselineskip
{\it Proof of Theorem \resIIa.} Assume $S_g$ is weakly compact for
every $g\in C(K).$ Suppose $E\subset X^*$ is a bounded subset which fails
to be \rwc. We must show there exists a \cseq\ which fails to tend
to zero uniformly on $E.$ Without loss of generality we may
assume $E=\seq{x^*}$ for some sequence $\seq{x^*}$ and there exists some
$\rho>0$ such that $E$ satisfies the conclusion of~\resIId\ with respect to
$\rho.$ \verify

Let $\mu_n\subset M(K)$ be a Hahn-Banach extension of $x^*_n$ and
let $\nu_n=|\mu_n|$. Let $\nu$ be a weak-star accumulation point
of $\brace{\nu_n}$ so that $\nu\geq0.$ Let $C=\sup\norm{\mu_n}.$
Choose $\delta_n>0$ so $\sum\delta_n<\rho/2.$ \verify

Let $U:X\to\lone(\nu)$ be the natural
inclusion. Then $U$ is weakly compact by the uniform integrability
criterion for weak compactness in $\lone(\nu)$. It now follows from
\resIId\ that there exists a sequence $\seq h$ in $X$ with $\norm{h_n}\leq1$
such that $\int|h_n|\,d\nu\goesto0$ and
$$\bigmod{x^*_{j_n}(h_n)}>\rho$$ for all $n$ and some sequence $\seq j.$
Choose $n_1$ so that $$\int|h_{n_1}|\,d\nu<{\delta_1\over2}.$$ We may
now find an increasing sequence $\pseq kl$ with $k_1=j_{n_1}$ so
that $$\int|h_{n_1}|\,d\nu_{k_l}<\delta_1\quad\hbox{for\ }l\geq2$$
and $$\left|\int h_{n_1}\,d\mu_{k_1}\right|>\rho.$$

Define $f_1=h_{n_1}.$ After renumbering we may now assume we have

\begingroup
\parindent=0pt
(1) $f_1\in B_X$.

(2) $\modulus{\int f_1\,d\mu_1}>\rho$.

(3) $\int |f_1|\,d|\mu_k|<\delta_1$ for $k>1.$
\endgroup

Let $g_1=1-|f_1|$ and redefine $\nu$ to be a weak-star
accumulation point of the new sequence $\seq\nu=\brace{|\mu_n|}$
which is now a subsequence of the sequence we started with.
Define $T:X\to C/X\llplus\lone(\nu)$
by $T=S_{g_1}\oplus U$ where $U$ is the operator defined
above. Since $S_{g_1}$ is weakly compact
by assumption it follows that $T$ is weakly compact. By \resIId\
there exists a new sequence $\seq h$ in $X$
with $\norm{h_n}\leq1$ such that $\norm{Th_n}\goesto0$ and
$$\bigmod{x^*_{j_n}(h_n)}>\rho$$ for all $n$ for some
sequence $\brace{j_n}.$ Note that these elements $\brace{x^*_{j_n}}$
are now being chosen from a subsequence of the original set $E$.
It is critical here that \resIId\  allows us to use the same constant
$\rho$ that
we used for the set $E.$

We now have
$$\dist\bigl(h_n(1-|f_1|),X\bigr)\goesto0$$ and
$$\int|h_n|\,d\nu\goesto0.$$ Choose $n_2$ so
that $\int|h_{n_2}|\,d\nu<\delta_2/2$ and
$$\dist\bigl(h_{n_2}(1-|f_1|\bigr),X)<{\delta_2\over2C}.$$ We may
now find an increasing sequence $\brace{k_l}$ with $k_0=1$ and $k_1=j_{n_2}$ so
that $$\int|h_{n_2}|\,d|\mu_{k_l}|<\delta_2\quad\hbox{for\ }l\geq2$$ and
$$\modulus{\int h_{n_2}\,d\mu_{k_1}}>\rho.$$

Define $f_2=h_{n_2}.$ After renumbering we may assume we have

\begingroup
\parindent=0pt
($1'$) $f_n\in B_X$ for $n=1,2.$

($2'$) $\modulus{\int f_n\,d\mu_n}>\rho$ for $n=1,2.$

($3'$) $\int|f_n|\,d|\mu_k|<\delta_n$ for $n=1,2$ and $k>n$.

($4'$) $\dist\bigl(f_2(1-|f_1|),X\bigr)<\delta_2/2C.$
\endgroup

Now let $g_2=(1-|f_1|)(1-|f_2|)$ and repeat the process.
At the $\Nth$ step we will have the following.

\begingroup
\parindent=0pt
($1''$) $f_n\in B_X$ for $1\leq n\leq N$.

($2''$) $\modulus{\int f_n\,d\mu_n}>\rho$ for $1\leq n\leq N$.

($3''$) $\int|f_n|\,d|\mu_k|<\delta_n$ for $k>n$ and $1\leq n\leq N$.

($4''$) $\dist\biggl({\wucprodii fj,X}\biggr)<\delta_n/2C$ for $2\leq n\leq N$.
\endgroup
\vskip\baselineskip

We now proceed as in \refBourPP, whose proof was elucidated in \refWojBS.
At the
$\Nth$ step define $\omega_N=\mu_N$. Let $\ph_1=f_1$ and $\ph_n=\wucprodii fj$
for $n>1$ so
$$|\ph_n|=\prod_{j=1}^{n-1}\parenth{1-|f_j|}-\prod_{j=1}^n\parenth{1-|f_j|}$$
and $\sum|\ph_n|\leq2$. Hence $\seq\ph$ is a w.u.C. series.

\begingroup %*%

\def\ssmprod{\prod_{j=1}^{n-1}\parenth{1-|f_j|}}
\def\om{\omega}

In general if $0\leq\alpha_j\leq 1$
then $1-\prod^s_{j=1}(1-\alpha_j)\leq\sum^s_{j=1}\alpha_j.$ Therefore for
$n\geq2$ we have
$$\eqalign{\modulus{\int \ph_n\,d\om_n}&=\modulus{\int\wucprodii fj\,d\om_n}\cr
&\geq\modulus{\int f_n\,d\om_n}-\modulus{\int
f_n\parenth{1-\ssmprod}\,d\om_n}\cr
&\geq\rho-\int\parenth{1-\ssmprod}\,d|\om_n|\cr
&\geq\rho-\int\sum^{n-1}_{j=1}|f_j|\,d|\om_n|\cr
&\geq\rho-\sum^{n-1}_{j=1}\delta_j\cr
&\geq{\rho\over2}.}$$

Choose~$\psi_n\in X$ with~$\psi_1=\ph_1$ and
$$\norm{\ph_n-\psi_n}\leq{\delta_n\over2C}$$ for~$n\geq 1.$
Then~$\sum\psi_n$ is a w.u.C. series and furthermore
$$\eqalign{
\modulus{\int\psi_n\,d\om_n}&\geq{\rho\over2}-{\delta_n\over2}\cr
&\geq{\rho\over4}}$$
for~$n\geq1.$ Since the sequence $\seq\om$ consists of Hahn-Banach extensions
of some sequence in $E$, it now follows from the notes at the beginning of
Section 1
that there exists a~\cseq\ in $X$ failing to tend to zero uniformly
on~$E$.\qed

\endgroup %*%
\vfil\eject

\ \section{3. Tight Uniform Algebras and Separable Distortions}

\comment{Transitory remarks. Remind them what the operator $S_g$ is.}

We will now discuss tightness and some of its connections to separably
distorted
dual spaces, Gleason parts, reducing bands and pointwise bounded approximation.

\proclaim{\resIIIa} Let $K$ be a compact space and let $X$ be a closed subspace
of $C(K)$. If~$\bd$ is a reducing band for~$X$ with~$\bd\subseteq\bxperp$ then
$${\bd\over\bd\cap\xperp}=\overline{\bigcup_{g\in C(K)}S^*_g(\bd\cap\xperp)}.$$
\endproclaim

\proof We claim that
$$\bd=\overline{\big\{g\,d\nu\,\bigm|\,\nu\in\bd\cap\xperp,\,g\in C(K)\big\}}.
\eqno\eqnIIIa$$ Let~$E$ be the right-hand side of \eqnIIIa\
so~$E\subseteq\bd.$ Let~$\mu\in\bd.$ Then, since~$\bd\subseteq\bxperp,$
it follows from a result of Chaumat (see Proposition V.17.11 in \refCon)
that there
exists some~$\nu\in\xperp$ such that~$\mu\ll\nu.$ Let~$\ldecomp\nu$ be the
\leb\
decomposition of~$\nu$ \wrt\ $\bd.$ Since~$\bd$ is a reducing band the
measure~$\nu_a$
lies in $\bd\cap\xperp.$ We now have~$\mu\ll\ldecomp\nu,\,\nu_a\perp\nu_s$
and~$\mu\perp\nu_s.$ Therefore~$\mu\ll\nu_a.$ \verify
Write~$d\mu=F\,d\nu_a$ for
some~$F\in\lone(\nu_a)$ and let~$\seq g$ be a sequence in~$C(K)$ so
that~$g_n\,d\nu_a\goesto d\mu$ in norm. Evidently~$\mu$ is in~$E$ which
implies~$\bd=E.$

Note that $S_g^*(\nu)=g\,d\nu+\xperp$ and
$S_g^*(\xperp)\subseteq\bmodxperp$. Since
$\bd$ is a reducing band the space $\bd/\bd\cap\xperp$ can be identified
isometrically with a closed subspace of $\bmodxperp$. The lemma now follows
from
Equation \eqnIIIa.
\qed

The following generalizes a result from \refMe.

\proclaim{\resIIIb} If $K$ is a metric space and $X$ is a tight subspace of
$C(K)$
then $\bmodxperp$ is separable.
\endproclaim

\proof Since $K$ is metrizable we may find a dense sequence $\seq g$ in $C(K)$.
\resIIIa\ now implies that
$$\bmodxperp=\overline{\bigcup_{n=1}^\infty S^*_{g_n}(\xperp)}.$$ Since weakly
compact sets in the dual of a separable Banach space are norm separable and
$S^*_g$ is weakly compact for all $g$, the result follows.\qed

Recall that if $A$ is a uniform algebra and $\net\bd$ is the collection of
all \ntmrb s
then $\bigllplus\bd_\alpha$ is isometric to a closed subspace of $A^*$ and
every nontrivial Gleason part corresponds to a distinct $\bd_\alpha$. We
therefore have the following result which is not new but was proved in
\refCG. However, the present proof is more elementary.

\proclaim{\resIIIr} If $A$ is a tight uniform algebra on a metric space $K$
then
$A$ has at most countably many nontrivial minimal reducing bands and at most
countably many nontrivial Gleason parts.
\endproclaim

The claim about the Gleason parts follows from the more basic fact that
$\norm{\varphi-\psi}=2$ for points $\varphi$ and $\psi$ in distinct parts.

For example if $A=A(\bidisk)$ is the bi-disk algebra, then
$\brace{z\times\Delta}$ is a nontrivial Gleason part for every $z$ on the unit
circle. In particular, the bi-disk algebra is not tight. For $R(K)$ where $K$
is a compact planar set, the fact about Gleason parts is well known. Any part
of $R(K)$ containing a non-peak point has positive area (see \refGamUA).

The only ingredient needed in the corollary is the separability of
$\bmodaperp$. We would like to mention that this is a special case of a more
general phenomenon. We say a \bs\ $Y$ is a {\it\sepdistone\ } if $Y=\mpl$ where
$M$ is separable and $L=\lone(\mu)$. Since every band is isomorphic to
$\lone(\mu)$ for some $\mu$, $A^*$ will be isomorphic to a \sepdistone\
whenever $\bmodaperp$ is separable. The following theorem now extends the
concept in the corollary.

\proclaim{\resIIIc} Let $A$ be a uniform algebra and suppose $A^*$ is
isomorphic to a closed subspace of a \sepdistone. Then $A$ has at most
countably
many non-trivial minimal reducing bands and therefore at most countably many
non-trivial minimal Gleason parts.
\endproclaim

\begingroup  % #1

\def\ea{{E_\alpha}}
\def\fa{{F_\alpha}}

\def\fnot{{F_{\alpha_\subnot}}}

\def\esum{{\bigllplus\kern-2pt\ea}}

\def\tinv{{\inverse T}}
\def\winv{{\inverse W}}
\def\uinvi{{U^{-1}_{\alpha_1}}}
\def\uinvii{{U^{-1}_{\alpha_2}}}
\def\qinv{{\inverse Q}}
\def\mpl{{M\llplus\kern-2pt L}}
\def\tint{{\norm T\norm\tinv}}
\def\d{{\hbox{d}}}
\def\ua{{U_\alpha}}
\def\uai{{U_{\alpha_1}}}
\def\uaii{{U_{\alpha_2}}}
\def\salph{{S_\alpha}}
\def\qm{{\hbox{$q_{\lower1.5pt\hbox{$\scriptscriptstyle M$}}$}}}
\def\ql{{\hbox{$q_{\lower1.5pt\hbox{$\scriptscriptstyle L$}}$}}}

This type of phenomenon has its origins in the paper \refHeni\ of
G.M. Henkin where it is shown that the ball-algebras~$A(\nball)$
are not isomorphic to the polydisk algebras~$A(\mpoly)$ when~$m$
is greater than one (also, see \refPeli). Our result is a direct extension
of Henkin's
work. We begin with some lemmas.

If~$\bd$ is any band then~$\hinf(\bd)$ is uniform algebra on its maximal ideal
space (see \hbox{Section 1).} It is easy to see that $\hinf(\bd)$ will be a
proper uniform algebra on its maximal ideal space if and only if it fails to be
self-adjoint. Note that if $\bd$ is a minimal reducing band then $\bd$ is
nontrivial if and only if $\bd\cap\aperp\ne0$. The following lemma shows that
when~$\bd$ is a minimal reducing band then~$\hinf(\bd)$ is a proper uniform
algebra on its maximal ideal space \iff\ $\bd$ is non-trivial.

\proclaim{\resIIId} Let~$A$ be a uniform algebra and let~$\bd$ be
a reducing band. Then \tfae

(a)~$\hinf(\bd)=\overline{\hinf(\bd)}.$

(b)~$\hinf(\bd)=\linf(\bd).$

(c)~$\bd\cap\aperp=0.$
\endproclaim

\proof

(a\iffarrow b) Assume~$\hinf(\bd)=\overline{\hinf(\bd)}.$ Note that
$$\hinf(\bd)=\left\{F\in\linf(\bd)\,\Bigg|\int F_\mu\,d\mu=0\hbox{\ for\ }
\mu\in\bd\cap\aperp\right\}.\eqno\eqnIIIb$$ Therefore, if~$f\in A$
then~$\int\bar f\,d\mu=0$ for~$\mu\in\bd\cap\aperp.$ Now, if~$g\in A$
and~$\mu\in\bd\cap\aperp$ then~$g\,d\mu\in\bd\cap\aperp.$ It then follows
that if~$\mu\in\bd\cap\aperp$ then~$\int\bar fg\,d\mu=0$ for all~$f,g\in A.$
Therefore~$\mu=0$ by the Stone-Weierstrass Theorem so~$\bd\cap\aperp=0.$ It
now follows from~\eqnIIIb\ that~$\hinf(\bd)=\linf(\bd).$

(b\iffarrow c) This follows immediately from \eqnIIIb.\verify\qed

The next result is a generalization of the fact that the
space~$L^1/H^1_0$ is not isomorphic to a subspace of
an \lonesp. The main ingredient is the theorem of Kisliakov from \refKisUA\
which
states that no proper uniform algebra is isomorphic to a quotient of a $C(K)$
space.

\proclaim{\resIIIg} Suppose~$\bd$ is a \ntmrb\ for some uniform
algebra~$A.$ Then~$\bmodbaperp$ is not isomorphic to a subspace of an \lonesp.
\endproclaim

\proof Suppose~$T:\bmodbaperp\to\lone(\mu)$ is an isomorphic embedding. Let~$E$
be a compact space such that~$C(E)$ is the dual of~$\lone(\mu).$
Then~$T^*:C(E)\to\hinf(\bd)$ is surjective. Since~$\bd$ is non-trivial we
have~$\bd\cap\aperp\ne0$. It then follows from \resIIId\
that~$\hinf(\bd)$ is a proper uniform algebra on its maximal ideal space which
contradicts Kisliakov's theorem.\qed

If~$X$ and~$Y$ are \bs s we say~$X$ is {\it~$C$-finitely
representable} in~$Y$ if for every finite-dimensional subspace~$F\subseteq X$
there exists a finite dimensional subspace~$G\subseteq Y$ such
that~$\bmd(F,G)\leq C$ where $d$ is the Banach-Mazur distance. If there exists
such a~$C$ we will simply say~$X$ is {\it finitely representable} in~$Y.$
We are
motivated by the following.

\proclaim{\resIIIe\ (Lindenstrauss-\pel, \refLP)} Suppose~$X$ is a \bs\
which is
finitely representable in~$\lone(\mu)$ for some $\mu$. Then~$X$ is
isomorphic to a
subspace of~$\lone(\mu')$ for some $\mu'$.
\endproclaim

We now study products which embed into separable distortions.

\proclaim{\resIIIf} Let~$\net E_{\alpha\in I}$ be a collection
of \bs s and let~$X=\esum.$ Suppose~$M$ is a separable \bs\ and
let~$L=\lone(\mu)$ be some~$\lone$-space. Assume there exists an isomorphic
embedding~$T:X\to\mpl$ and let~$C=\tint.$ If~$I_\subnot$ is the set
of~$\alpha$ in~$I$ such that~$\ea$ fails to be~$2C$-\fr\ in~$L$
then~$I_\subnot$ is countable.
\endproclaim

\proof Assume~$I_\subnot$ is uncountable. We may then assume that~$I$ is
uncountable and we have an isomorphic embedding~$T:X\to\mpl$
where~$E_\alpha$ fails to be~$2C$-\fr\ in~$L$ for every~$\alpha\in I$
where~$C=\tint.$ Therefore, for every~$\alpha\in I$ we may find a finite
dimensional subspace~$\fa\subseteq\ea$ such that~$$\d(\fa,G)>2C\eqno\eqnIIIc$$
for every subspace~$G$ of~$L$ such that~$\dim G=\dim \fa.$

Since~$I$ is uncountable we may assume \Wlog\ that there exists a
fixed integer~$n$ independent of~$\alpha$ such that~$\dim \fa=n$
for all~$\alpha\in I.$ Choose~$\eps>0$ so
that~$${2\over 2-\eps}(1+2\eps)={3\over 2}.\eqno\eqnIIId$$
It is well-known that the Banach-Mazur distance
on the space of~$n$-dimensional \bs s is a separable metric. We may therefore
assume \Wlog\ that~$$\d(\fa,F_{\alpha'})\leq1+\eps\eqno\eqnIIIe$$ for
all~$\alpha$ and~$\alpha'$ in~$I.$

Let~$\alpha_\subnot$ be any element of~$I.$ For every~$\alpha\in I$ let
$\ua:\fnot\biggoesto\fa$ be an isomorphism with
$$\norm\ua\norm{\inverse\ua}\leq1+2\eps,\eqno\eqnIIIf$$ which can
be done by \eqnIIIe. Furthermore, after multiplying by a constant we
may assume~$$\norm{\inverse\ua}=\eps\eqno\eqnIIIg$$ for all~$\alpha.$

Let~$\qm$ and~$\ql$ be the natural projections from $\mpl$ to $M$ and $L$,
respectively.
%\ifnum\showdiagrams=1
%$$\diagram
%&\mpl\dlto_\qm\drto^\ql&\\
%M& &L.
%\enddiagram
%$$
%\else\nodiag
%\fi
For every~$\alpha\in I$ let~$i_\alpha:\fa\injects X$ be the natural
injection and define~$S_\alpha:\fnot\biggoesto M$ by
$\salph=\qm\circ T\circ i_\alpha\circ\ua.$ Note that the space of
bounded linear operators~$L(\fnot,M)$ is
separable. Since~$\net S_{\alpha\in I}$ is
an uncountable collection in~$L(\fnot,M)$ we may find distinct
elements~$\alpha_1$ and~$\alpha_2$ in~$I$ so that
$$\norm{S_{\alpha_1}-S_{\alpha_2}}<{1\over\norm\tinv}.
\eqno\eqnIIIh$$ Define~$W:\fnot\biggoesto X$ by
$W=(i_{\alpha}\circ U_{\alpha_1}-i_{\alpha}\circ U_{\alpha_2}).$

\pclaim1 We have~$$\norm W\norm\winv\leq1+2\eps.\eqno\eqnIIIi$$
If~$x\in\fnot$ then
$$\eqalignno{\norm{Wx}&=\norm\uai+\norm\uaii\cr
 &\geq\norm x\left({1\over\norm\uinvi}+{1\over\norm\uinvii}\right)\cr
 \noalign{\hbox{so}}
 \norm\winv&\leq{\norm\uinvi\norm\uinvii\over\norm\uinvi+\norm\uinvii}\cr
 \noalign{\vskip4pt}
 &={\eps\over2}&\eqnIIIj}$$
by \eqnIIIg. Since~$\norm W\leq\norm\uai+\norm\uaii$ we have
$$\eqalign{\norm W\norm\winv&\leq\norm\uaii\norm\uinvii\left(
 {\norm\uinvi\over\norm\uinvi+\norm\uinvii}\right)\cr
 &\quad+\norm\uai\norm\uinvi\left(
 {\norm\uinvii\over\norm\uinvi+\norm\uinvii}\right)\cr
 &\leq1+2\eps}$$
by \eqnIIIf.

Define~$Q:\fnot\biggoesto L$ by~$Q=\ql\circ T\circ W.$

\pclaim2 We have
$$\norm Q\norm\qinv\leq{3\over2}C.\eqno\eqnIIIk$$
If~$x\in\fnot$ then
$$\eqalignno{\norm{Qx}&=\norm{\ql(TWx)}\cr
 &=\norm{(TW)(x)}-\norm{(\qm TW)(x)}\cr
 &=\norm{(TW)(x)}-\norm{(S_{\alpha_1}-S_{\alpha_2})(x)}\cr
 &\geq{1\over\norm\tinv\norm\winv}\norm x-{1\over\norm\tinv}\norm x\cr
 &={\norm x\over\norm\tinv}\left({1\over\norm\winv}-1\right)\cr
 \noalign{\hbox{so}}
 \norm\qinv&\leq{\norm\tinv\norm\winv\over1-\norm\winv}.}$$
Furthermore, since~$\norm Q\leq\norm T\norm W$ we have
$$\eqalignno{\norm Q\norm\qinv&\leq C{\norm W\norm\winv\over1-\norm\winv}\cr
 &\leq{2\over2-\eps}C\norm W\norm\winv\cr
 &\leq{2\over2-\eps}(1+2\eps)C\cr
 &={3\over 2}C}$$
by \eqnIIId.

If we let~$G=Q(\fnot)$ then \eqnIIIk\ implies~$\d(\fnot,G)\leq{3\over2}C.$
This contradicts \eqnIIIc. Hence, \resIIIf\ is proved.\qed

\endgroup % #1

By \resIIIe\ we have the following corollary.

\proclaim{\resIIIt} If~$\net E$ is a collection of \bs s such that
the product~$\bigllplus E_\alpha$ embeds isomorphically into a \sepdistone\
then all but a countable number of the~$E_\alpha$ embed isomorphically
into some~$\lone$-space (where the~$\lone$-space depends on~$\alpha$).
\endproclaim

We now summarize our results in the following proof.

{\it Proof of \resIIIc.} Assume~$A$ is a uniform algebra and~$A^*$
is isomorphic to subspace of an \sepdistone. If~$\net\bd$ is
the collection of all \ntmrb s then the sum~$\mrbsum$ is isometric to a
subspace
of~$A^*.$ By \resIIIt\ this
implies~$\bd_\alpha/\bd_\alpha\cap\aperp$ embeds in some \lonesp\ for
all but a countable number of~$\alpha.$ However, \resIIIg\ states
that~$\bd_\alpha/\bd_\alpha\cap\aperp$ fails to embed in an \lonesp\ for
every~$\alpha\in I.$ Therefore, the set~$I$ must be countable.
Furthermore, every non-trivial Gleason part corresponds to a distinct
\ntmrb, which finishes the proof of the proposition.\qed

At the present time it is not known if there exists a uniform algebra $A$
with the
property that $\bmodaperp$ is nonseparable and the dual of $A$ embeds into
a separable
distortion of an $\lone$-space. One problem is that it is not clear when
$\bd/\bd\cap\aperp$ is separable for an arbitrary \ntmrb\ $\bd$. Even if
this problem
is solved, the \ntmrb s do not exhaust $\bmodaperp$. A discussion of the
complete
decomposition of $\baperp$ can be found in \refCG.

The separability of $\bmodaperp$ can be applied to some ideas in
pointwise approximation. To illustrate, let~$K$ be a compact planar set and
let~$R(K)$ be the space of functions in~$C(K)$ which can be uniformly
approximated by rational functions with poles off~$K.$ Define~$Q\subset K$ to
be the collection of non-peak points for~$R(K)$ and let~$\lamq$ be the
restriction of planar \leb\ measure to~$Q.$ Define~$\hinf(\lamq)$ be the
weak-star
closure of~$R(K)$ in~$\linf(\lamq).$ It is a theorem of A.M. Davie in
\refDav\ that
if~$f\in\hinf(\lamq)$ then there exists a sequence of functions~$\seq{f}$
in~$R(K)$
such that~$f_n\goesto f$ \pwae~[$\lamq$] and~$\norm{f_n}\leq\norm f.$
This conclusion is sometimes referred to as pointwise bounded
approximation with a reduction in norm. It is known that Davie's theorem
implies, without much difficulty, that every $z\in Q$ has a representing
measure absolutely continuous with respect to $\lamq$.

We take the following approach to this problem (also, see \refCG\ or \refCon).
Let $A$ be an arbitrary uniform algebra. Given~$m\in\baperp$ we have the
natural
projection~$$\Thbatoh(m)\eqno\eqnIIIl$$ which is the dual of the
injection~$$\Svarlaperptoba m\eqno\eqnIIIm$$ (see Section 1). Because the
space~$\hinf(\baperp)$ is identified isometrically with a closed subspace
of~$\ddual A$ it follows from Goldstine's Theorem that if~$F\in\hinf(\baperp)$
with~$\norm F\leq1$ then there exists a net~$\net f$ in~$A$
with~$\norm\falpha\leq1$
such that~$\falpha\goesstar F.$ Let~$m$ be any measure in~$\baperp.$ Recall
that
$\bmodaperp$ contains all the point evaluations for the points off the Choquet
boundary. Therefore this net also has the property that~$\falpha(z)\goesto
f(z)$ for
every~$z\in Q$ where $Q$ is the complement of the Choquet boundary. Since the
natural projection of~$\hinf(\baperp)$ into~$\hinf(m)$ is weak-star
continuous, we
have~$\falpha\goesstar F$ in~$\hinf(m)$ (where we are using~$F$ as the
symbol for an
element of~$\hinf(\baperp)$ as well as its projection~$F_m$
into~$\hinf(m)).$ It is
now easy to see that there exists a sequence~$\seq f$ bounded by one such
that~$f_n(z)\goesto F(z)$ for all~$z\in Q$ and~$f_n\goesto F$
\pwae~[$m$]. We therefore have the following.

\proclaim{\resIIIh} If~$A$ is a uniform algebra and~$m\in\baperp$
then for any~$F$ in the unit ball of~$\hinf(\baperp)$ there exists a
sequence~$\seq f$ from the unit ball of~$A$ such that~$f_n\goesstar F$
in~$\hinf(\baperp)$ and~$f_n\goesto F$ \pwae~[m].
\endproclaim

Recall that a linear operator $T$ between the \bs s~$X$
and~$Y$ is a quotient map if the induced injection~$S:X/Z\to Y,$
where~$Z=\ker T,$ is an isometry. The following proposition relates the Davie
phenomenon directly to the projection $\tau$.

\proclaim{\resIIIi} Let~$A$ be a uniform algebra
on an arbitrary compact space~$K$ and let~$m\in\baperp.$ Then \tfae

(a) For every~$f\in\hinf(m)$ there exists a sequence~$\seq f$
in~$A$ with~$\norm{f_n}\leq\norm f$ such that~$f_n\goesto f$ \pwae~[m].

(b) The natural projection~$\Thbatoh(m)$ is a quotient map.
\endproclaim

\comment{can we replace $m$ by $m_a$ in part (b), and have $m$ be arbitrary?}

\proof (a\implies b) Assume that (a) holds and let~$I$ be the
kernel of~$\tau$ and~$$S:\hinf(\baperp)/I\biggoesto\hinf(m)$$ be the induced
injection. Given~$f\in\hinf(m)$ let~$\seq f\subset A$ be
the sequence mentioned in the statement of~(a). Let~$F$ be a
weak-star accumulation point of~$\seq f$
in~$\hinf(\baperp)$ so~$\norm F\leq\norm f.$ Since the map~$\tau$ is a dual
map it is continuous from the weak-star topology to the weak-star
topology. Therefore,~$\tau(f_n)$ accumulates weak-star at~$\tau(F)$
in~$\hinf(m)$ and~$\seq f$ converges
weak-star to~$f$ so~$\tau(F)=f.$ Hence the map~$S$ is onto.
Furthermore,we have~$\norm{\tau(F)}\leq\norm F$ by definition
so~$\norm f=\norm F.$ Hence~$S$ is an isometry and~$\tau$ is a
quotient map and therefore (b) holds.

(b\implies a) Assume~$\tau$ is a quotient map and
let~$f\in\hinf(m).$ We may then find an~$\widetilde F\in\hinf(\baperp)$
so that~$\tau(\widetilde F)=f$ and~$\norm{\widetilde F+I}=\norm f$
where~$I=\ker\tau.$ Since~$I$ is a weak-star closed subspace
of~$\hinf(\baperp)$ we may find an~$F\in\hinf(\baperp)$ such
that~$\tau(F)=\tau(\widetilde F)=f$ and $\norm{\widetilde F+I}=\norm F.$
Hence,~$\norm F=\norm f.$ By \resIIIh\ we may find
a sequence~$\seq f$ in~$A$ such that~$\norm{f_n}\leq\norm F$
and~$f_n$ converges to~$f$ \pwae~[$m$] which is the desired conclusion.\qed

It is possible that~$\tau$ is a quotient map \iff\ it is onto, but this is
currently not known to be true. If ~$\tau$ were onto then note
the kernel of~$\tau$, call it~$I$, is an ideal and the induced
map~$\hinf(\baperp)/I\,\buildrel\tilde\tau\over\longrightarrow\hinf(m)$
is an algebra isomorphism from a Banach algebra to a uniform algebra. If
$\hinf(\baperp)/I$ is a uniform algebra then~$\tilde\tau$ would be an
isometry and so~$\tau$ would be a quotient map. However,
it is not clear when~$\hinf(\baperp)/I$ is a uniform algebra.

Nevertheless, the conclusion of Davie's theorem can now be stated in
terms of specific properties of the natural projection. If $A$ is a uniform
algebra on a compact space $K$ and $m\in\baperp$ we say $m$ is an {\it
ordinary Davie measure} if $\tau$ is a quotient map and $m$ is a {\it strong
Davie measure} if $\tau$ is an isometry. In general, a linear operator
between \bs
s is an isometric embedding \iff\ its dual is a quotient map. Therefore,
$m$ is an
ordinary Davie measure \iff\ $\sigma$ is an isometric embedding and is
strong Davie
measure \iff\ $\sigma$ is an isometry. Since $\tau$ is an algebra homomorphism
between uniform algebras, $\tau$ will be an isometry as soon as it is an
isomorphism. Since $\sigma$ is always an injection, we see that $m$ is a strong
Davie measure as soon as $\sigma$ is onto. Since the evaluations for the
points off
the Choquet boundary lie in $\bmodaperp$, it follows easily that when $m$
is a strong
Davie measure then every point of the Choquet boundary has a representing
measure
absolutely continuous \wrt\ $m$.

For the sake of completeness we will briefly discuss the injectivity of
$\tau$. We say $m\in M(K)$ is a {\it weakly rich (resp. strongly rich) measure}
for $A$ if when $\seq f$ is a bounded sequence in $A$ such that
$\int|f_n|\,d|m|\goesto0$ then $f_ng+A\goesweak0$ (resp.
$\norm{f_ng+A}\goesto0$) for every $g\in C(K)$. The concept of richness was
introduced in \refBourDPP\ where it was shown that a uniform algebra $A$
(or even an
arbitrary subspace of $C(K)$) has the \pp\ if there exists a
strongly rich measure for $A$. Note that weakly rich measures on strongly tight
spaces (where the operators $S_g$ are compact) are strongly rich.

The following result can be found in \refMe. If $m\in M(K)$ let $m=m_a+m_s$ be
the \leb\ decomposition of $m$ \wrt\ $\baperp$.

\proclaim{\resIIIj}  Let~$A$ be a uniform algebra on a
compact space~$K$ and let~$m$ be an element of~$M(K)$. Then \tfae

(a) The natural projection~$\Thbatoh(m_a)$ is one-to-one.

(b) If~$\{f_n\}$ is a bounded sequence in~$A$ such
that~$\int|f_n|\,d|m|\goesto0$ then~$f_n\goesstar0$ in~$\linf(\mu)$ for
every~$\mu\in\aperp$.

(c)~$m$ is a weakly rich measure for~$A$.

Furthermore, if the above hold, and~$K$ is metrizable, then~$\bmodaperp$ is
separable.
\endproclaim

We will now show that $A$ possesses a strong Davie measure whenever
$\bmodaperp$ is separable. We approach this problem from a general point of
view. Let~$X$ be a \bs\ and let~$E\subseteq X$ be a closed subspace. If~$Y$
is a
subspace of~$X$ we say~$Y$ is a {\it full subspace \wrt\ $E$} if the induced
map~$$\YtoE\eqno\eqnIIIn$$ is an isometric embedding; that is, for every~$y\in
Y$ we have $$\dist(y,Y\cap E)=\dist(y,E).\eqno\eqnIIIo$$

Note that we do not assume~$Y$ to be a closed subspace of~$X$ which means
that~$Y/Y\cap E$ may only have a semi-norm. Therefore, when we say
isometric embedding in the definition, what we really mean is that~$\sigma$
preserves the semi-norm.

The case we should be thinking about is ~$X=\baperp,$ the space~$E$
is~$\aperp,$ and~$Y=\lone(m)$ for some measure~$m\in\baperp.$ When~$\lone(m)$
is full with respect to~$\aperp,$ then~$m$ is an ordinary Davie
measure. If the map~$\sigma$ is onto,~$m$ will be a strong Davie measure.
When the subspace~$E$ is clear we will simply refer to~$Y$ as being
a full subspace. Furthermore, we will identify~$Y/Y\cap E$ with its
image~$\sigma(Y/Y\cap E)$ in~$X/E.$ Note that all of our \bs s are complex and
$\clspan F$ refers to the closed complex linear span of $F$.

\proclaim{\resIIIl} Let~$X$ be a \bs, and let~$E$ be a closed
subspace of~$X.$ Suppose~$S$ is a separable subset of~$X$ (respectively,
of~$X/E$). Then there exists a closed, full, separable subspace~$Y\subseteq X$
such that~$S\subseteq  Y$ (respectively,~$S\subseteq Y/Y\cap E$).
\endproclaim

{% begin1

\def\normYE{{\norm{x_n+Y\cap E}}}
\def\normE{{\norm{x_n+E}}}
\def\YE{{Y\cap E}}
\def\YmodE{Y/\YE}

\proof Let~$\seq s$ (respectively~$\brace{s_n+E}$) be a dense sequence
in~$S$ and let
$${\seq x}^\infty_{n=1}=\left\{\sum\limits^N_{k=1}\alpha_ks_k\,\Biggl |\,
 \alpha_k=p+iq,\quad p,q\in{\Bbb Q},\quad N\in{\Bbb N}\right\}.$$
Choose~$x_{n,k}\in E$ so that
$\lim\limits_{k\to\infty}\norm{x_n+x_{n,k}}=\normE.$ Let
$Y=\clspan\bigl\{\seq x\cup\brace{x_{n,k}}\bigr\}$ so that~$S\subseteq Y$
(respectively,~$S\subseteq\YmodE$). Now, since~$x_{n,k}\in\YE,$
$$\normYE\leq\lim\limits_{k\to\infty}\norm{x_n+x_{n,k}}=\normE\leq\normYE$$
so~$\normYE=\normE$ for all~$n$. Hence, the map~$\sigma$ in
Equation~\eqnIIIn\ is an isometry on~$\brace{x_n+\YE}.$ By definition, the
sequence~$\seq x$ is dense in~$\clspan\seq x$ and therefore
$\brace{x_n+\YE}$ is
dense in~$\YmodE.$ Hence $\sigma$ is an isometric embedding since it
preserves the
norm on a dense set.
\qed
}% end1

If~$X$ is a band of measures and~$m\in X$ we identify~$\lone(m)$ with the
subband of~$X$ consisting of all measures absolutely continuous \wrt\ to~$m.$

\proclaim{\resIIIm} Let~$X$ be a band of measures on some
compact metric space~$K$ and let~$E\subseteq X$ be a closed subspace.
Suppose~$S\subseteq X/E$ is separable. Then there exists an~$m\in X$ such
that~$\lone(m)$ is full and~$S\subseteq\lone(m)/\lone(m)\cap E.$
\endproclaim

\proof By \resIIIl\ we can find a closed, full, separable
subspace~$Y_1\subseteq X$ with~$S\subseteq Y_1/Y_1\cap E.$ Let~$\bd_1$ be the
band generated by~$Y_1$ so that~$\bd_1$ is separable (here we use the
metrizability of $K$). Using ~\resIIIl\ again we may find a closed,
full, separable subspace~$Y_2$ with~$\bd_1\subseteq Y_2.$ Let~$\bd_2$ be the
band generated by~$Y_2$ and repeat, so we have
$$Y_1\subseteq\bd_1\subseteq Y_2\subseteq\bd_2\subseteq\cdots$$
where~$Y_n$ is a closed, full, separable subspace and~$\bd_n$ is a separable
band.

{%begin stuff
\def\ub{{\overline{\bigcup\limits_{n=1}^\infty \bd_n}}}
\def\YE{Y\cap E}

Let~$$\bd=\ub$$ so~$\bd$ is separable. It is easy to see $\bd$ is a band.

\claim~$\bd$ is full. Note that~$\bd=\overline{\bigcup Y_n}.$
Let~$Y=\bigcup Y_n.$ If~$y\in Y$ then~$y\in Y_n$ for some~$n$ and
$$\dist(y,\YE)\leq\dist(y,Y_n\cap E)=\dist(y,E)\leq\dist(y,\YE)$$
since~$Y_n$ is full. This shows~$Y$ is full which clearly implies~$\overline Y$
is full. Hence,~$\bd$ is full.

Since~$\bd$ is separable it follows that $\bd=\lone(m)$ for some $m\in X$ and
the proposition is proved.\qed

}%end stuff

\proclaim{\resIIIn} Suppose~$A$ is a uniform algebra on a compact
metric space~$K$ and~$G$ is any subset of~$K$ which does not meet the Choquet
boundary. Then \tfae

(a)~$G$ is separable in the Gleason metric.

(b) There exists an ordinary Davie measure~$m$ such that every
point in~$G$ has a representing measure absolutely continuous with respect
to~$m.$
\endproclaim

\proof Let~$G'$ be the subset of~$A^*$ consisting of point
evaluations at the points in~$G$ so that $G'\subseteq\bmodaperp.$

(a\implies b) Assume (a) holds. Since the Gleason metric on~$G$
corresponds to the norm on~$\bmodaperp$ it follows that ~$G'$ is a separable
subset of~$\bmodaperp.$ Using the above proposition with~$X=\baperp$ and
$E=\aperp$ we may find a measure~$m$ in~$\baperp$ such that~$\lone(m)$
is a full subspace and~$$G'\subseteq\varlaperp m.$$
Since~$\lone(m)$ is full~$m$ is an ordinary Davie measure.
Furthermore, if~$z\in G$ then there exists a representing measure~$\mu$
for~$z$ such that~$\mu+\aperp\in G'.$ We may then find a~$g\in\lone(m)$ such
that $\mu+\aperp=g\,dm+\aperp$ and so~$g\,dm$ is a
complex representing measure for~$z$. We can then find a representing
measure for~$z$ absolutely continuous \wrt\ $m.$ This proves (b).

(b\implies a) follows from the fact then~$\lone(m)$ is separable when $K$
is metrizable (we don't need $m$ to be an ordinary Davie measure here).\qed

The next corollary is immediate.

\proclaim{\resIIIo} If~$A$ is a uniform algebra on a compact
metric space and~$\bmodaperp$ is separable then~$A$ admits a strong Davie
measure~$m.$
\endproclaim

Applying \resIIIj\ we have another corollary.

\proclaim{\resIIIq} Let~$A$ be a uniform algebra on a compact
metric space~$K.$ Then \tfae

(a)~$A$ admits a weakly rich measure.

(b)~$A$ admits a strong Davie measure.

(c)~$\bmodaperp$ is separable.
\endproclaim

>From \resIIIb\ we deduce the following.

\proclaim{\resIIIp} If $A$ is a tight uniform algebra on compact metric
space $K$ then $A$ admits a strong Davie measure $m$. In particular, every
point
off the Choquet boundary has a representing measure absolutely
continuous \wrt\ $m.$
\endproclaim

Obviously $\lambda_Q$ is a strong Davie measure for $R(K)$. It follows from
\resIVb\ that if $D$ is a smoothly bounded strictly
pseudoconvex domain in $\Cn$ then the surface-area measure on $\boundaryD$
is also
a strong Davie measure for $A(\boundaryD)$. It is known that $R(K)$ has the
\pp\ for
any compact planar set $K$ and $A(\boundaryD)$ has the \pp\ for the domains
$D$ just
mentioned (see the comments in Sections 1 and 5). It follows that these uniform
algebras have \wsc\ duals. If $A^*$ is \wsc\ then $\bmodaperp$ is \wsc. If
$m$ is an
ordinary Davie measure then $\lone(m)/\lone(m)\cap\aperp$ is isomorphic to
a subspace
of $\bmodaperp$ and is therefore \wsc\ and has $\hinf(m)$ as its dual. This
proves
the following theorem, which has an analogous version for $R(K)$ and is a
direct
generalization of a theorem in \refMoon\ which deals with the unit disk.

\proclaim{\resIIIv} Suppose $D\subset\Cn$ is a bounded strictly
pseudoconvex domain
with $C^2$ boundary. Let $m$ be surface-area measure on $\partial D$.
Suppose $\seq
f$ is a sequence of functions in $\lone(m)$ such that $\lim\int f_nh\,dm$
exists for
every $h\in\hinf(m)$. Then there exists an $f\in\lone(m)$ such that
$\lim\int f_nh\,dm=\int fh\,dm$ for every $h\in\hinf(m)$.
\endproclaim

\section{4. Uniformly Convergent Fourier Series}

Let $\Gamma$ be the unit circle in $\C$ and let $U$ be the space of continuous
functions $F$ on $\Gamma$ which extend to be analytic in the unit disk such
that
the series $\sum_{n=0}^\infty\hat F(n)z^n$ converges uniformly to $F$ where
$\hat
F(n)={1\over2\pi}\int^\pi_{-\pi} F(\exp{i\theta})\exp{-in\theta}\,d\theta$. For
$g\in L^1(d\theta)$ let $P_n(g)=\sum_{k=0}^n\hat g(k)z^k$. If we define
$$\norm{F}_U=\sup_{n\geq0}\norm{P_n(F)}_\infty$$ then $U$ becomes a Banach
space
with the norm $\norm{\cdot}_U$. In this section we will prove the following.

\proclaim{\resIVna} Let $U$ be the space of analytic uniformly convergent
Fourier
series on the unit circle with the above norm. Then $U$ has the following
properties.

(a) $U$ is not isomorphic to a quotient of $C(G)$ for any compact space $G$.

(b) $U$ has the \pp.

(c) $U^*$ is weakly sequentially complete.

(d) $U$ and $U^*$ have the \dpp.

(e) $U^*$ is isomorphic to a separable distortion of an $\lone$-space.
\endproclaim

\varproof{Proof of \resIVna\ (a)} Suppose, on the contrary, there exists a
compact
space $G$ and a surjective continuous linear operator $T:C(G)\to U$. Let
$A$ be the
disk algebra on $\Gamma$ and let $t:U\to A$ be the natural inclusion. Let
$P:A\to
l^2$ be the Paley operator $P(f)=(\hat f(2^n)_{n=0}^\infty$. Then $P$ is
1-summing
(see \refPeli). However $PtT$ is also 1-summing and therefore compact
(again, see
\refPeli), which implies by the surjectivity of $T$ that $Pt$ is compact. By
examining $Pt(z^{2^n})$, we see that this is a contradiction.
\qed

It is well-known that part (c) of \resIVna\ follows from (b). To prove (b)
we will show
that $U$ embeds isometrically into some $C(K)$ space as a tight subspace
and apply
\resIIa.

Let $K'=\brace{1/n}_{n=1}^\infty\cup\brace{0}$ and let $K=K'\times\Gamma$.
Define a
sequence of closed subspaces of $K$ by
$\Gamma_n=\brace{1/(n+1)}\times\Gamma$ for
$n\geq0$ and let $\Gamma_\infty=\brace{0}\times\Gamma$. If $\Phi\in C(K)$
we can write
$\Phi=(\vf_\infty,\vf_0,\vf_1,\vf_2,\ldots)$ where
$\restrict\Phi;{\Gamma_n}=\vf_n$ for $n\geq0$,
$\restrict\Phi;{\Gamma_\infty}=\vf_\infty$ and
$\vf_n\goesto\vf_\infty$ uniformly on $\Gamma$. Define an isometry $i:U\to
C(K)$ by
$i(F)=(F,P_0(F),P_1(F),P_2(F),\ldots)$ and let $X=i(U)$.

To show $X$ is a tight subspace of $C(K)$ we must study the operators $S_g$
defined
in Section 1. Given a compact space $G$ and closed subspace $Y\subseteq
C(G)$, we
define $Y_{cg}$ to be the set of $g\in C(K)$ such that $S_g$ is weakly
compact (the
``c" represents Brian Cole and the ``g" represents Ted Gamelin). It was
shown in
\refCG\ that when $A$ is a uniform algebra on $G$ then $A_{cg}$ is a closed
subalgebra of $C(G)$. However, the proof does not use the algebraic
structure of
$A$ and goes as follows. From general theory we know that $S_g:Y\to C/Y$ is
weakly
compact \iff\ $\ddual{S_g}(\ddual Y)\subseteq C/Y$. Since
$\ddual{S_g}(F)=gF+\ddual
Y$, $S_g$ is weakly compact \iff\ $g\ddual Y\subseteq\ddual Y+C$. It is now
clear that $Y_{cg}$ is a closed subalgebra of $C(G)$.

Let $Y_{CG}$ be the set of those $g$ such that $S_g$ is compact. It is an even
easier task to show that $Y_{CG}$ is again a closed subalgebra of $C(G)$.
This result
is proved in \refMe\ for algebras, however the proof does not use the algebraic
structure.

The fact that $X_{cg}$ is a closed subalgebra of $C(K)$ means that we need
only verify $S_g$ is weakly compact on a set of continuous functions $g$
that is
self-adjoint and separates the points of $K$.

\proclaim{\resIVnb} $X$ is a tight subspace of $C(K)$.
\endproclaim

\proof
By the observation above, we need only verify that $S_\Phi$ is weakly
compact for
$\Phi=(\vf_\infty,\vf_0,\vf_1,\vf_2,\ldots)$ where $\Phi$ has the following
form.
First, suppose $\vf_n=\delta_{nm}z^{\pm1}$ for some integer $m$ and
$n\geq0$ where
$\delta_{nm}$ is the Kronecker delta function. Then
$S_\Phi$ is easily seen to be of finite rank. Secondly, suppose $\vf_n=z$ for
$n\geq0$ or $\vf_n=\bar z$ for $n\geq0$. If we prove $S_\Phi$ is weakly
compact in
this case then, since the functions in the first and second cases form a
separating
self-adjoint family, we will have shown $X_{cg}=C(K)$ by the comments
above. This,
by definition, means $X$ is a tight subspace of $C(K)$.

Suppose $\vf_n=z$ for $n\geq0$. Suppose $f\in X$ and $f=i(F)$ where $F\in
U$. Then if
$a_n=\hat F(n)$ we have
$$\eqalign{
\Phi f &= (zF,a_0z,a_0z+a_1z^2,\cdots)\cr
&=(zF,0,a_0z,a_0z+a_1z^2,\cdots)+(0,a_0z,a_1z^2,a_2z^3,\cdots).}
$$
Therefore, if we
define $V:X\to C(K)$ by $Vf=(0,a_0z,a_1z^2,a_2z^3,\dots)$ then $V$ is
easily seen
to be continuous and $S_\Phi=qV$ where $q:C(K)\to C(K)/X$ is the natural
quotient
map. Furthermore, if we let $j:X\to l^2$ map $f$ to its Fourier
coefficients and
let $\tilde V:l^2\to C(K)$ by $\tilde Vx=(0,x(0)z,x(1)z^2,\ldots)$, then
$j$ and
$\tilde V$ are continuous and $V=\tilde Vj$. Hence $V$ is weakly compact which
implies $S_\Phi$ is weakly compact. The argument for $\vf_n=\bar z$ is
similar.\qed

\varproof{Proof of \resIVna\ (b), (c) and (e)} In \resIIa\ we proved that tight
subspaces have the \pp, and therefore so does $X$. Part (b) now follows
from the
fact that $U$ is isomorphic to $X$. It is well-known that Banach spaces with
the \pp\ have weakly sequentially complete dual spaces, which takes care of
(c).
Part (e) is an immediate consequence of \resIIIb.\qed

If $G$ is any compact space and $Y\subseteq C(G)$ is a closed subspace, let
$Y_b$
and $Y_B$ be the space of functions $g\in C(G)$ such that $S_g$ (respectfully,
$S_g^{**}$) is completely continuous. These are called the {\it Bourgain
algebras}
of $Y$. These spaces were first defined in \refCT. It is not difficult to
see that
$Y_b$ and $Y_B$ are closed subalgebras of $C(G)$, as was shown in \refCT. The
motivation to study these spaces was the work of Bourgain in \refBourDPP.
It can be
deduced from Bourgain's work that if $Y_B=C(K)$ then $Y$ and $Y^*$ have the
\dpp. This
is how we plan to prove part (d) of \resIVna.

The lemma below follows immediately from \resIIc.

\proclaim{\resIVnc} Suppose $A$, $B$ and $C$ are \bs s with continuous linear
operators $S:A\to B$ and $T:A\to C$. Assume $T$ is weakly compact and
$\ddual T$ is
completely continuous. Suppose $T$ has the additional property that
whenever $x_n$
is a bounded sequence in $A$ such that $\norm{Tx_n}\goesto0$ then $\norm
{Sx_n}\goesto0$. Then $\ddual S$ is completely continuous.
\endproclaim

\proclaim{\resIVnd} The Bourgain algebra $X_B$ of $X$ equals $C(K)$.
\endproclaim

\proof Since $X_B$ is a closed subalgebra of $C(K)$, it suffices to show
$\ddual{S_\Phi}$ is completely continuous for the family of functions $\Phi$
studied in the proof of \resIVnb. We need only consider the functions
$\Phi=(z,z,z,\ldots)$ and $\Phi=(\bar z,\bar z,\bar z,\ldots)$.

Let $m=\restrict{{d\theta\over2\pi}};\Gamma_\infty$; that is, $m\in M(K)$ is
normalized \leb\ measure on the set $\Gamma_\infty$. Let $T':C(K)\to
L^1(m)$ be the
natural inclusion. Then $T'$ is weakly compact. Since $C(K)^{**}$ has the \dpp,
$(T')^{**}$ is completely continuous. Let $T$ be the restriction of $T'$ to
$X$. Let
$V$ be the operator defined in the proof of \resIVnb. Then
$$\norm{Vf}\leq{1\over2\pi}\int |f|\,dm.\eqno{\eqnIVna}$$ Since
$S_\Phi=qV$, we see that $S_\Phi$ and $T$ satisfy the hypothesis of
\resIVnc. Hence,
$\ddual{S_\Phi}$ is completely continuous. As in \resIVnb, the proof for
$\Phi=(\bar
z,\bar z,\bar z,\ldots)$ is similar.\qed

A comment is in order. The inequality \eqnIVna\ implies that the
operator $V$ is 1-summing (strictly 1-integral, in fact; see \refWojBS\ for the
definitions). Therefore $S_\Phi$ is 1-summing. This provides us with
another way of
deducing the above properties of $S_\Phi$. It is well-known that the second
adjoint of a 1-summing operator is 1-summing, and that 1-summing operators
are weakly
compact and completely continuous.

\varproof{Proof of \resIVna\ (d)} In \resIVnd\ we showed that
$\ddual{S_\Phi}$ is
completely continuous for every $\Phi\in C(K)$. It can now be deduced from
the work in
\refBourDPP\ that $X$ and $X^*$ have the \dpp. Since $U$ is isomorphic to
$X$, the
proof is finished.\qed

The space $X$ is tight but it is not strongly tight; that is, the operators
$S_\Phi$ are not compact for every $\Phi\in C(K)$. We will make this result
precise in the proposition below. This is interesting because in every
known example
where $A$ is tight uniform algebra on some compact space $G$, $A$ turns out
to be
strongly tight. It if not known if this is true in general.

\proclaim{\resIVne} The operator $S_\Phi$ is compact \iff\
$\restrict\Phi;{\Gamma_\infty}$ is constant. That is,
$$X_{CG}=\{\Phi\in C(K) \,\big |\, \restrict\Phi;{\Gamma_\infty}\hbox{\rm\ is
constant}\}.$$

\endproclaim

\proof If $Y$ is the right-hand side of the above, is not hard to see that
$Y\subseteq X_{CG}$. If $\vf_n=\delta_{nm}z^{\pm1}$ as in the beginning of
the proof of
\resIVnb\ then
$S_\Phi$ is compact. Since this is true for a self-adjoint family of
functions which
separates the points of
${\seq\Gamma}_{n=0}^\infty$, and $X_{CG}$ is a closed subalgebra of $C(K)$,
$X_{CG}$ contains all functions which are constant on $\Gamma_\infty$.

Let $\phi_n=i(z^n)$ so that ${\seq\phi}_{n=0}^\infty$ is a sequence in the
unit ball of $X$. Let $\Phi\in C(K)$ be arbitrary. \claimII
$S_\Phi(\phi_n)\goesweak0$. Since $S_\Phi$ is weakly compact it suffices to
show
that zero is the only weak accumulation point of $\brace{S_\Phi(\phi_n)}$.
So we
let $\Psi+X$ be such a weak accumulation point. If we write
$\Phi=(\vf_\infty,\vf_0,\vf_1,\vf_2,\ldots)$ then
$\Phi\phi_n=(z^n\vf_\infty,0,0,\ldots,0,z^n\vf_n,z^n\vf_{n+1},\ldots)$. Write
$\Psi=(\psi_\infty,\psi_0,\psi_1,\ldots)$. If $\nu_0$ is any measure on
$\Gamma$
with $\int\,d\nu_0=0$ then $\restrict\nu_0;{\Gamma_0}\in X^\perp$ which implies
$\int\psi_0\,d\nu_0=0$ and hence $\psi_0=c_0$ is constant. Similarly we can
show
$\psi_1=c_{1,1}+c_{1,2}z$. If
$\mu=\restrict{{d\theta\over2\pi}};{\Gamma_0}-
\restrict{{d\theta\over2\pi}};{\Gamma_1}$
then $\mu\in X^\perp$ and $\int_K\Psi\,d\mu=0$ which implies $c_0=c_{1,1}$.
Similarly, by considering the annihilating measure
$\restrict{\bar z{d\theta\over2\pi}};{\Gamma_1}-
\restrict{\bar z{d\theta\over2\pi}};{\Gamma_2}$, we see that
$\psi_1=c_0+c_1z$ and
$\psi_2=c_0+c_1z+c_2z^2$. We repeat and find that $\Psi$ is in $X$, which
proves
the claim.

Now, suppose $\Phi\in X_{CG}$. Let $N$ be a positive integer. Let
$\nu_n=\restrict{{{\bar z^{n+N}\over2\pi}d\theta}};\Gamma_n$ so that
$\norm{\nu_n}=1$ and $\nu_n\in X^\perp$. Then $\int_K
\Phi\phi_n\,d\nu_n=\widehat{\vf_n}(N)$ which implies
$$|\widehat{\vf_n}(N)|\leq\norm{\Phi\phi_n+X}.$$ Since $S_{\Phi}(\phi_n)$ tends
to zero weakly and $S_{\Phi}$ is compact it follows that
$\norm{S_\Phi(\phi_n)}\goesto0$ and so $\widehat{\vf_\infty}(N)=0$. If $N$ is a
negative integer then we prove $\widehat{\vf_\infty}(N)=0$ by using the
annihilating measure
$\restrict{\bar z^{n-|N|}{d\theta\over2\pi}};{\Gamma_{n-1}}-
\restrict{\bar z^{n-|N|}{d\theta\over2\pi}};{\Gamma_n}$ for $n\geq1$. Hence,
$\vf_\infty$ is constant and we are done.\qed

\section{5. A Note on Inner Functions}

We conclude with an application of tightness to inner functions on strictly
pseudoconvex domains.

Let $D$ be a domain in~$\Cn$ with smooth boundary and let~$A=A(\boundaryD)$
be the
uniform algebra of functions in~$C(\boundaryD)$ which extend to be analytic
in~$D.$  Let~$m$ be the normalized surface-area measure on~$\boundaryD$ and
let~$\hinf(\boundaryD)=\hinf(m)$ be the corresponding Hardy space. For short
we will write~$C=C(\boundaryD)$ and~$\hinf=\hinf(\boundaryD).$ Recall that a
function~$f\in\hinf(\boundaryD)$ is an inner function if~$|f|=1$ a.e. [$m$].

\proclaim{\resIVa} Let $D$ be a strictly pseudoconvex domain with
$C^2$ boundary in $\Cn$. Suppose~$f$ is an inner function
in~$\hinf(\boundaryD).$ If $f(z_n)\goesto0$ for some sequence $\seq z$ tending
to $\boundaryD$ then $\barf\notin\hinf+C$. In particular, if $n>1$ then
$\barf\in\hinf+C$ \iff\ $f$ is constant. If $D$ is the unit disk then
$\barf\in\hinf+C$ \iff\ $f$ is a finite Blaschke product.
\endproclaim

The proof is indirect and utilizes the \pp\ and the theory of tight uniform
algebras. As discussed in Section 1, it was proven in \refCG\ (also, see
\refMe) that if $D$ is a strictly pseudoconvex domain with $C^2$ boundary
in $\Cn$
then $A(D)$ is a strongly tight uniform algebra. Actually all that is
needed is the
solvability of the $\dbar$-problem with H\"older estimates on the
solutions, and
therefore this is all that is needed in \resIVa. It is proven in \refMe\ that
whenever a uniform algebra is strongly tight on some compact space $K$, it is
strongly tight as a uniform algebra on its Shilov boundary. Hence,
$A(\boundaryD)$ is strongly tight on $\boundaryD$. Since strongly tight uniform
algebras are tight, it follows from \resIIa\ that $A(\boundaryD)$ has the
\pp\ (this
result was also in \refMe).

We will need to lift the properties of the operators $S_g$ on
$A(\boundaryD)$ to
the corresponding operators on the uniform algebra $\hinf$. We accomplish this
by the following result, which can be found in \refCG\ (this can also be
deduced
from the results in \refCR).

\proclaim{\resIVb} The measure~$m$ is in~$\baperp$ and the natural
projection~$\Thbatoh(m)$ is an isometry.
\endproclaim

Consider~$\hinf(m)$ as a uniform algebra on the maximal ideal
space of~$\linf(m).$ Given a function~$g\in\linf(m)$
define~$$S_{g,\hinf}:\hinf(m)\biggoesto {\linf(m)\over\hinf(m)}$$ by~$$f\mapsto
fg+\hinf(m).$$ We define~$(\hinf(m))_{CG}$ to be set of those~$g$ such
that~$S_{g,\hinf}$ is compact.

\proclaim{\resIVc}~$\hinf(m)+C(\boundaryD)\subseteq(\hinf(m))_{CG}.$
\endproclaim

\proof Let~$g\in C(\boundaryD).$ Then~$S_g:A\to C/A$ is compact.
Now~$$S_g^*:\aperp\biggoesto\varBaplusSa$$
satisfies~$S_g^*(\aperp)\subseteq\bmodaperp$ and is given by~$\nu\mapsto
g\,d\nu+\aperp.$ Let ~$p:A^*\to\bmodaperp$ be the natural map so
that~$p\circ S^*_g$ is compact. Let
$$T_g=\bigrestrict{(p\circ S^*_g)};{\lone(m)\cap\aperp}$$ so
$$T_g:\lone(m)\cap\aperp\biggoesto\varbmodaperp$$ is compact.

By \resIVb\ the natural projection~$\tau$ is an isometry and
therefore its predual~$$\Svarlaperptoba m$$ is also an isometry.
Let~$U_g=\sigma^{-1}\circ T_g$ so~$U_g$ is compact and
$$U_g:\lone(m)\cap\aperp\biggoesto\varlaperp m$$
by~$$f\,dm\mapsto fg\,dm+\lone(m)\cap\aperp.$$ The
adjoint of~$U_g$ is $S_{g,\hinf}$ which is therefore compact.
Hence,~$g\in(\hinf)_{CG}$ and we have shown
that $C(\boundaryD)\subseteq(\hinf)_{CG}.$ Since~$S_{g,\hinf}=0$
for~$g\in\hinf$ the
proposition is proved.\qed

\varproof{Proof of \resIVa} Assume~$f$ is an inner function
and~$f(z_n)\goesto0$ for some sequence $\seq z$ in $D$ such that $z_n\goesto
z$ where $z\in\boundaryD$. We will show~$\barf\not\in(\hinf(m))_{CG}$; in
other words,~$S_{\barf,\hinf}$ fails to be compact. It then follows from
\resIVc\ that~$\barf\not\in\hinf+C.$ Our method will be to show that
$S_{\barf,\hinf}$ is an isomorphism on a copy of $c_0$ in $\hinf$.

\begingroup % #1

\def\phi{{\varphi}}

Let $\seq\ph\subset A^*$ be the evaluation functionals at $\seq z$. Since it is
well-known that every point on $\boundaryD$ is a peak point for $A$ it
follows that
$\seq\phi$ is not relatively weakly compact. Since $A(\boundaryD)$ has the
\pp\ we
may find, after passing to a subsequence of $\seq z$ if necessary, a \wuc\
series~$\Sigma\psi_n$ in~$A$ and an $\eps>0$ such that
$|\psi_n(z_n)|>\eps$ for all $n.$ We claim that
$\norm{\barf\psi_n+\hinf}\notgoesto0.$ Otherwise there exist bounded
sequences~$\seq j\subset \hinf$ and~$\seq k\subset\linf(m)$
with~$\norm{k_n}\goesto0$
such that~$\barf\psi_n+j_n=k_n.$ Therefore, since $f$ is inner,
$\psi_n+fj_n$ tends
to zero uniformly which contradicts that fact that~$f(z_n)\rightarrow0$.

It is now clear that $S_{\barf,\hinf}$ fails to be compact. In fact, it
follows from the remarks in Section 1 that $S_{\barf,\hinf}$ is an isomorphism
on a copy of $c_0$.\qed

\endgroup % #1

{\bf Acknowledgments.} Many of the ideas in this paper were conceived
while the author was a graduate student at Brown University under the
supervision of
Brian Cole. The author would also like to thank T.W. Gamelin for some helpful
discussions,  and thank the referee for demonstrating that the theory
of tight subspaces could be applied to the space of uniformly convergent
Taylor series.
\vskip1truein

\vfil\eject\

\def\ref#1;#2;#3;#4;#5;#6;{{\bibrmfont \item{#1} #2, {\bibitfont #3,}
#4 {\bibbffont #5} #6.}}
     %number,author(s),title,journal,volume,whatever.

\def\book#1;#2;#3;#4;{{\bibrmfont \item{#1} #2, {\bibitfont #3,} #4.}}
     %number,author(s),title,whatever.
\def\bysame{\hbox to 30pt{\hrulefill}}

\frenchspacing
\centerline{\csc References}
\vskip0.75in

\parindent=0pt

\ref\refBP;C. Bessaga and A. Pe\l czy\'nski;On bases and unconditional
convergence of series in Banach spaces;Studia. Math;17;(1958), pp. 151--16;

\book\refBourU;J. Bourgain; Quelques proprietes lineaires topologiques de
l'espace des
series de Fourier uniformement convergentes;Seminaire Initiation a
l'Analyse, G.
Choquet, M. Rogalski, J. Saint Raymond, 22e annee. Univ. Paris-6, 1982-83,
Expose no.
14;

\ref\refBourPP;\bysame;On weak completeness of
the dual of spaces of analytic and smooth
functions;Bull. Soc. Math. Belg. Serie B;35; (1) (1983), 111--118;

\ref\refBourDPP;\bysame;The \dpp\ for the
ball-algebras, the polydisc algebras and the Sobolev
spaces;Studia Math.;77;(3) (1984), 245--253;

%\ref\refBD;J. Bourgain and F. Delbaen;A class of special $\cal L_\infty$
%spaces;Acta. Math;145;(1980), 155--176;

%\ref\refCCG;T.K. Carne, B.J. Cole and T.W. Gamelin;A uniform algebra
%of analytic functions on a Banach space;Trans. Amer. Math. Soc.;314;(2)
%(1989), 639--659;

\ref\refCT;J.A. Cima and R.M. Timoney;The \dpp\ for certain planar
uniform algebras;Mich. Math J.;34;(1987), 99--104;

\ref\refCG;B.J. Cole and T.W. Gamelin;Tight uniform algebras and algebras
of analytic functions;J. Funct. Anal.;46;(1982), 158--220;

\ref\refCR;B.J. Cole and R.M. Range;A-measures on complex manifolds and
some applications;J. Funct. Anal.;11;(1972), 393--400;

\book\refCon;J.B. Conway;The Theory of Subnormal Operators;The American
Mathematical Society, Providence, 1991;

\ref\refDav;A.M. Davie;Bounded limits of analytic functions;Proc. Amer.
Math. Soc.;32;(1972), 127--133;

\ref\refDeli;F. Delbaen;Weakly compact operators on the disc algebra;J.
of Algebra;45;(1977), 284--294;

\ref\refDelii;\bysame;The \pel\ property for some uniform algebras;Studia
Math;64;(1979), 117--125;

%\book\refDS;N. Dunford and J.T. Schwartz;Linear Operators;Part I,
%Interscience, New York, 1958;

\book\refGamUA;T.W. Gamelin;Uniform Algebras;Chelsea Publishing Co., New
York, 1984;

\book\refHeni;G.M. Henkin;The Banach spaces of analytic functions in a sphere
and in a bicylinder are not isomorphic;Functional Anal. Appl. {\bf 2} (4)
(1968),
334--341\semicolon\ Funct. Anal. Prilo\u z. {\bf 2} (4) (1968), 82--91;

\ref\refJames;R.C. James;Bases and reflexivity of Banach spaces;Annals of
Math;52;(3) (1950), 518--527;

%\ref\refPrieto;J.A. Jaramillo and A. Prieto;Weak-polynomial convergence on
a Banach
%space;Proc. Amer. Math Soc.;118;(2) (1993), 463--468;

\item{\refKisPP} S.V. Kisliakov, {\it On the conditions
of Dunford-Pettis, \pel, and Grothendieck}, Soviet Math. Dokl.
{\bf 16} (1975), 1616--1621\semicolon\
Dokl. Akad. Nauk. SSSR {\bf 225} (1975), 1252--1255.

\ref\refKisUA;\bysame;Proper uniform algebras are uncomplemented;Dokl.
Akad. Nauk.-SSSR (Russian);309;(4) (1989), 795--798;

\ref\refLP;J. Lindenstrauss and A. \pel;Absolutely summing operators
in $\lone$-spaces and their applications;Studia Math.;29;(1968), 275--326;

\ref\refMilne;H. Milne;Banach space properties of uniform algebras;Bull.
London Math. Soc.;4; (1972), 323--326;

\ref\refMoon;M.C. Mooney;A theorem on bounded analytic functions;Pac. J.
Math.;43;(1973), 457--463;

\ref\refPP;A. \pel;Banach spaces on which every unconditionally converging
operator
is weakly compact;Bull. Acad. Polon. Sci.;10;(1962), 641--648;

\ref\refPeli;\bysame;Banach spaces of analytic functions and absolutely
summing operators;CBMS Regional
Conference Series in Mathematics, The American Mathematical Society,
Providence;30;(1977);

\ref\refMe;S.F. Saccone;Banach space properties of strongly tight uniform
algebras;Studia Math.;114;(2) (1995), 159--180;

\ref\refWojPP;P. Wojtaszczyk;On weakly compact operators from some uniform
algebras;Studia Math.;64;(1979),
105--116;

\book\refWojBS;\bysame;Banach Spaces for Analysts;Cambridge University Press,
New York, 1991;

\bye